\documentclass[11pt]{article}
\usepackage{graphicx,epsfig}
\usepackage{amsmath,amssymb,euscript}
\usepackage{latexsym}
\usepackage{color}

\newtheorem{theorem}{\bf Theorem}[section]

\newtheorem{lemma}[theorem]{\bf Lemma}
\newtheorem{proposition}[theorem]{\bf Proposition}

\newtheorem{es}{Example}[section]
\def\NN{{\mathbb N}}
\def\RR{{\mathbb R}}

\def\PP{{\mathbb P}}
\def\C{{\mathcal{C}}}

\newcommand{\cond}{\rm{cond}}
\begin{document}
\title{Filtered interpolation for solving \\ Prandtl's  integro-differential equations}
\author{M.C. De Bonis, D. Occorsio, W. Themistoclakis {\thanks{corresponding author}}}

\maketitle
\begin{abstract}
In order to solve Prandtl--type equations we propose a collocation--quadrature method based on VP filtered interpolation at Chebyshev nodes. Uniform convergence and stability are proved in a couple of H\"older--Zygmund spaces of locally continuous functions. With respect to classical methods based on Lagrange interpolation at the same collocation nodes, we succeed in reproducing the optimal convergence rates of the $L^2$ case by cutting off the typical $log$ factor which seemed inevitable dealing with uniform norms. Such an improvement does not require a greater computational effort. In particular we propose a fast algorithm based on the solution of a simple 2-bandwidth linear system and prove that, as its dimension tends to infinity, the sequence of the condition numbers (in any natural matrix norm) tends to a finite limit.
\end{abstract}



\section{Introduction}
 In this paper we propose a numerical method for  the following Prandtl-type equation
\begin{eqnarray}\nonumber \sigma f(y)&-&\frac 1 \pi \int_{-1}^1\frac{f(x)}{(x-y)^2}\varphi(x)dx-\frac 1 \pi \int_{-1}^1 \log |x-y| f(x)\varphi(x) dx\\ &+&\frac 1 \pi \int_{-1}^1 h(x,y)f(x)\varphi(x)dx= g(y),\quad y\in (-1,1),\label{iniziale2}\end{eqnarray}
where $\varphi(x)=\sqrt{1-x^2}$, $\sigma$ is a real constant, $h$ is a smooth kernel function and the first integral has to be understood as the following derivative
\begin{equation}\label{der}
\int_{-1}^1\frac{f(x)}{(x-y)^2}\varphi(x)dx=\frac{d}{dy} \int_{-1}^1\frac{f(x)}{x-y}\varphi(x)dx,
\end{equation}
being the integral at the right--hand side a Cauchy principal value integral.

 Integro Differential Equations (IDEs) of the type (\ref{iniziale2}) are models  for many physics and engineering problems (see \cite{vainikko} and the references therein). Indeed, besides the well known Prandtl's Equation which  governs the circulation air flow along the contour of a plane wing profile (see e.g. \cite{dragos1,dragos2}), Prandtl-type equations (\ref{iniziale2}) model, for instance,   the load transfer problem from  thin-walled elements to massive bodies (see e.g. \cite{mkhi,koer}),  water scattering problems of vertical or inclined plates merged in infinitely deep water (see e.g. \cite{Paul}),  crack problems in  composite material
or in  inhomogeneous bodies \cite{tvad}, etc.


We focus on the case that the right--hand side known--term $g$ can be well approximated by polynomials. More precisely, for a given Jacobi weight $u$, we consider the error of best, Jacobi--weighted, uniform approximation of $g$ by means of polynomials of degree at most $n$, namely
\[
E_n(g)_u=\inf_{\deg (P)\le n}\|(g-P)u\|, \qquad   \|g u\|:=\sup_{|x|\le 1}\left|g(x)u(x)\right|,
\]
and we suppose that
\begin{equation}\label{hp-g}
E_n(g)_u=O\left(\frac 1{n^s}\right), \qquad \mbox{ for a certain $s>0$}.
\end{equation}
The space of all functions satisfying (\ref{hp-g}) is here denoted by $Z_s(u)$ and it is a Banach space equipped with the norm
\[
\|g\|_{Z_s(u)}:=\|gu\|+\sup_{n> 0}(1+n)^s E_n(g)_u,\qquad s>0.
\]
Such kind of B--spaces have been introduced in \cite{DT_Besov, JuLu-JCAM97}.  They are equivalent to  some H\"older--Zygmund--type spaces defined by means of Ditzian-Totik weighted $\varphi$--moduli of smoothness \cite{DT, DT_Besov} and several singular integral equations have been already studied in a couple of these, also known as Besov--type spaces (see e.g. \cite{JuLu-JCAM97, CaCrJuLu00, debonis2014, mastroDB2006, DBO, operatoreD, MaRuThLAG, air})

We prove that for any $g\in Z_s(\varphi)$, the equation (\ref{iniziale2}) has a unique solution $f\in Z_{s+1}(\varphi)$ (cf. Theorem \ref{uni}). Such a result can be proved by standard arguments, combining the Fredholm's alternative theorem with the analysis of the mapping properties of the operators involved in the equation \cite{CaCrJuLu00, DBO}.

The main result of this paper is the construction of a numerical method in order to approximate the solution $f$ by a polynomial $f_n$  such that

\[
\|(f-f_n)\varphi\|\le \frac\C{n^{s}}\|g\|_{Z_{s}(\varphi)}
\]
and
\[
\|(f-f_n)\varphi\|_{Z_r(\varphi)}\le \frac \C{n^{s-r}}\|g\|_{Z_{s}(\varphi)}, \qquad 0<r\le s,
\]
where, in both the estimates, $\C$ denotes a positive constant independent of $n$ and $f$.

The proposed numerical method falls into the class of polynomial projection methods (see e.g. \cite{Ak,PS} ), but it is based on  non--standard projections obtained by applying de la Vall\'ee Poussin (briefly VP) filters. They project the function into non--standard polynomial spaces, spanned by the so--called fundamental VP polynomials that share the interpolation property of fundamental Lagrange polynomial at Chebyshev nodes of second type \cite{Ca.Th-wave, ThL1, woula_NA, woula_99, ThBa}. The resulting filtered VP interpolation has been already applied to the field of Cauchy singular integral equation in order to solve the airfoil equation \cite{air} and in other situations \cite{mastroDB2006, operatoreD}. It results to be an useful device also in our case. The related approximation error has been recently characterized in \cite{filtered1D} and error estimates in H\"older--Zygmund spaces can be found in \cite{papero-mata}.

Here, we are going to use the already known results on the near--best uniform approximation  provided by filtered VP interpolation, in order to construct a convergent, stable and efficient method for the numerical solution of (\ref{iniziale2}). In the particular case that $h(x,y)\equiv 0$, the algorithm fast computes the numerical solution $f_n$ by solving a linear system with a band matrix whose not null $(i,j)$--entries are those such that $|i-j|\in\{0,2\}$ (cf. (\ref{sist-H0})). Moreover, we prove that, as $n\rightarrow \infty$, the condition numbers of such  systems converge to a finite limit, whatever natural matrix norm we consider (cf. Theorem \ref{th-kappa}).

Several authors studied Prandtl's type equations by proposing different numerical approaches. Among them we recall \cite{calio, CaCrJu97, CaCrJuLu00, DBO}. In \cite{calio} a Nystr\"om method based on the local approximation by quasi-interpolatory splines is proposed for the particular case $h\equiv 0$. Nevertheless, the efficiency of  this method seems to suffer from  the low degree of approximation provided by splines functions. A global polynomial approximation based on the Lagrange interpolation at Jacobi zeros has been considered in \cite{CaCrJu97, CaCrJuLu00, DBO}. In \cite{CaCrJu97} a collocation and a collocation--quadrature method  is studied in Sobolev subspaces of weighted--$L^2$ spaces, getting optimal convergence rates. For the same methods, convergence results in  uniform norms have been obtained in \cite{CaCrJuLu00} by investigating a regularized version of the equation in certain weighted Besov  spaces of locally continuous functions.
A direct collocation--quadrature method has been proposed in \cite{DBO}, considering the equation in the H\"older--Zygmund spaces $Z_s(u)$. In such spaces the authors prove the convergence and stability of their method and  study also the conditioning of the final linear system. Nevertheless, dealing with uniform norms, the previous methods based on Lagrange interpolation suffer from the unboundedness of the Lebesgue constants, which entails the presence of  a $log$ factor in the corresponding error estimates. This extra--factor disappears in the method here proposed, thanks to the uniform boundedness of the VP interpolating operators w.r.t. weighted uniform norms.

The paper is organized as follows. Section 2 is a survey of the theoretical results on the Prandtl equation (\ref{iniziale2}): we analyze the spaces where the equation is studied and the related mapping properties of the continuous operators that compose the equation, arriving to state existence and uniqueness theorems for the solution of (\ref{iniziale2}). In Section 3 we focus on the VP interpolating polynomials we use for the numerical method, recalling the related properties of interest. The method is described in the next two sections. Section 4 is devoted to the approximate equation and its unique solvability, providing convergence and stability theorems. Section 5 deals with the computation of the numerical solution by analyzing the computational aspects of the linear system obtained by our method. Finally, Section 6 provides the numerical experiments, by comparing the results of our method with those achieved by the methods in \cite{calio, CaCrJuLu00, DBO}. For a better legibility, the proofs of the main results are given in the Appendix.

\section{The equation in H\"older-Zygmund spaces}
Defining the operators
\begin{eqnarray}
\label{D}
D:f\rightarrow Df, && Df(y):=-\frac 1 \pi \int_{-1}^1\frac{f(x)}{(x-y)^2}\varphi(x)dx,\\
\label{K}
K:f\rightarrow Kf, && Kf(y):=-\frac 1 \pi \int_{-1}^1 \log |x-y| f(x)\varphi(x) dx,\\
\label{H}
H:f\rightarrow Hf, && Hf(y):=\frac 1 \pi \int_{-1}^1 h(x,y)f(x)\varphi(x)dx,
\end{eqnarray}
and set $I f= f$, the equation (\ref{iniziale2}) can be shortly written as follows
\begin{equation}\label{eq-cont-op}
(\sigma I+ D+K+H)f=g.
\end{equation}
Here we are going to study this equation  and the related operators in some subspaces of the space $C^0_\varphi$ of all locally continuous functions $f$ (i.e. continuous in any $[a,b]\subset]-1,1[$) such that
\[
\lim_{x\rightarrow \pm 1} f(x)\varphi(x)=0.
\]
It is well--known that $C^0_\varphi$ is a Banach space equipped with the norm
$\|f\|_{C^0_\varphi}:= \|f\varphi\|$.
Moreover, Weierstrass approximation theorem holds in $C^0_\varphi$, i.e. any $f\in C^0_\varphi$ can be approximated by polynomials with the desired accuracy and we have
\[
f\in C^0_\varphi \Longleftrightarrow \lim_{n\rightarrow \infty}E_n(f)_\varphi=0.
\]
The rate of convergence to zero of $E_n(f)_\varphi$ can be characterized  by means of certain weighted moduli of smoothness that Ditzian and Totik introduced in \cite{DT} in order to get Jackson and Stechkin type inequalities.

In particular, for any $r\in\NN$, if we consider the following Sobolev--type  spaces
\[W_r(\varphi)=\left\{f\in C^0_\varphi: \mbox{$f^{(r-1)}$ is locally absolutely continuous and } \|f^{(r)}\varphi^{r+1}\|<\infty\right\},\]
equipped with the  norm
\[\|f\|_{W_r(\varphi)}=\|f\varphi\| + \|f^{(r)}\varphi^{r+1}\|,\]
then, taking into account the equivalence of Ditzian--Totik moduli with some $K$--functional \cite[Th. 2.1.1]{DT}, we get
\begin{equation}\label{Favard-Sob}
E_n(f)_\varphi\leq \frac{\C}{n^r} \|f\|_{W_r(\varphi)}, \quad \forall f\in W_r(\varphi), \quad r\in\NN,\qquad \C\neq \C(n,f),
\end{equation}
where, throughout the paper,  $\C$ denotes a positive constant having different meaning in different formulas and we  write $\C \neq \mathcal{C}(n,f,\ldots)$  to say that $\C$ is  independent of  $n,f,\ldots$.

In order to generalize (\ref{Favard-Sob}) to the case of any, not necessarily integer, exponent $r>0$, we introduce the following spaces
\[
Z_r (\varphi):=\{f\in C^0_\varphi \ :\ \sup_{n>0}(n+1)^rE_n(f)_\varphi<\infty\}, \qquad r>0,
\]
equipped with the following norm
\[
\|f\|_{Z_r(\varphi)}:=\|f\varphi\|+\sup_{n>0}(n+1)^rE_n(f)_\varphi, \qquad r>0.
\]
Obviously, by this definition, we have
\begin{equation}
\label{Favard-Zig}E_n(f)_\varphi\leq \frac{\C}{n^r} \|f\|_{Z_r(\varphi)},\qquad \forall f\in Z_r(\varphi),\quad r>0, \qquad \C\neq \C(n,f).
\end{equation}
Regarding the definition of the spaces $Z_r(\varphi)$, we recall they are also known as H\"older--Zygmund spaces and constitute a particular case of the Besov-type spaces studied in \cite{DT_Besov}, where they have been equivalently defined in terms of Ditzian--Totik moduli of smoothness (\cite[Th. 2.1]{DT_Besov}).

Moreover, we recall that the spaces $Z_r(\varphi)$ belong to the larger class of  spaces $C_\varphi^{\mathcal B}$ introduced in \cite[p. 204]{JuLu-JCAM97}, where a more general behavior of the error $E_n(f)_\varphi$ is allowed by an arbitrary sequence ${\cal B}:=\{b_n\}_n\subset \RR^+$ such that
\[
\lim_{n\rightarrow\infty}b_n=0\qquad \mbox{and}\qquad \sup_{n>0}\frac{E_n(f)_\varphi}{b_n}<\infty.
\]
We will consider also the non--weighted space $C^0=C[-1,1]$ equipped with the uniform norm, and its Zygmund--type subspaces $Z_r$, $r>0$, which are defined as $Z_r(\varphi)$ but replacing the function $\varphi$ with the constant $1$.

In \cite{JuLu-JCAM97} the following embedding properties have been proved
\begin{itemize}
\item [(i)] $\forall r>0$, $Z_r(\varphi)$ is compactly embedded in $C^0_\varphi$  cf. \cite[Lemma 3.2]{JuLu-JCAM97});
\item [(ii)] $\forall s>r>0$, $Z_s(\varphi)$ is compactly embedded in $Z_r(\varphi)$ (cf. \cite[Lemma 3.3]{JuLu-JCAM97});
\item [(iii)] $\forall s>0$, $Z_s$ is continuously embedded in $Z_s(\varphi)$ (cf. \cite[Remark 3.5]{JuLu-JCAM97}).
\end{itemize}
Moreover, recalling that the non weighted  error $E_n(f):=\inf_{\deg(P)\le n}\|f-P\|$ satisfies the Favard inequality (see e.g. \cite{Mhaskar-book})
\[
E_n(f)\le \frac\C n E_{n-1}(f^\prime)_\varphi,\qquad \C\ne\C(n,f),
\]
we easily deduce that
\begin{equation}\label{f-f'}
\sup_{n>0}(n+1)^{r+1}E_n(f)\le \C \|f^\prime\|_{Z_{r}(\varphi)}\le \C \|f^\prime\|_{Z_{r}}, \qquad  r>0, \qquad \C\ne\C(f).
\end{equation}
In the following subsections we are going to recall the main properties of the operarors $D,K$ and $H$ in the previous spaces in order to get, finally, the unique solvability of the Prandtl equation (\ref{eq-cont-op}) in H\"older--Zygmund spaces.
\subsection{\underline{On the operator $D$}}
Recalling (\ref{der}), the hypersingular integral operator $D$  is strictly related to the Cauchy singular integral operator
\begin{equation}\label{A}
Af(y):=-\frac 1 \pi \int_{-1}^1\frac{f(x)}{x-y}\ \varphi(x)dx
\end{equation}
being 
\begin{equation}\label{DA}
Df(y)=\frac d{dy}Af(y), \qquad -1<y<1.
\end{equation}
Consequently, some nice properties of the operator $A$ (see e.g. \cite{PS}) have been transferred on $D$.

In particular, let us consider the Chebyshev polynomials of the second kind of degree $n\in \NN$
  \begin{equation}\label{p}
 p_n(x)=\sqrt{\frac 2\pi}\ \frac{\sin[(n+1)t]}{\sin t}, \qquad t=\arccos x,\quad |x|\le 1,
 \end{equation}
being understood that $p_n(\pm 1)=\sqrt{\frac 2\pi}(n+1)$. In \cite{CaCrJu97} the following mapping property has been proved
\begin{equation}\label{inva-D}
D p_n(y)=(n+1) p_n(y), \qquad \forall n\in\NN.
\end{equation}
Moreover, in the recent paper \cite{DBO}, it has been proved the following
\begin{theorem}\label{th-D}
For any $r>0$,  $D:Z_{r+1}(\varphi)\rightarrow Z_r(\varphi)$ is a bounded  map having bounded inverse.
\end{theorem}
Finally, the explicit form of $D^{-1}: Z_r(\varphi)\rightarrow Z_{r+1}(\varphi)$ has been found in \cite[Prop. 2.3]{CaCrJuLu00}
\begin{equation}\label{D-1}
D^{-1} f(y)=-\hat{A}WAf(y), \qquad \forall f\in Z_r(\varphi),
\end{equation}
where $A$ is given by (\ref{A}) and
\begin{eqnarray}\label{A-1}
\hat{A}f(y)&:=&\frac 1 \pi \int_{-1}^1\frac{f(x)}{x-y}\frac{dx}{\varphi(x)},\\
\label{W}
Wf(y)&:=&\frac 1 \pi \int_{-1}^1\log|x-y|f(x)\frac{dx}{\varphi(x)}.
\end{eqnarray}
\subsection{\underline{ On the operator $K$}}
It is well known (see e.g. \cite[p. 8]{Ak}) that
\[
\sup_{|x|\le 1}\int_{-1}^1\log |x-y|dx<\infty.
\]
This yields
\begin{equation}\label{K-inf}
\|Kf\|\le\C\|f\varphi\|,\quad \forall f\in C^0_\varphi,\qquad \C\ne\C(f),
\end{equation}
namely, $K: C^0_\varphi\rightarrow C^0$ is a bounded map.

On the other hand, $K$ is related to the Cauchy operator $A$ given in (\ref{A}) as follows
\begin{equation}\label{KA}
Af(y)=\frac d{dy}Kf(y), \qquad -1<y<1.
\end{equation}
In \cite{operatoreD, mastroDB2006} the mapping properties of Cauchy singular integral operators have been studied. In particular, we recall that \cite[Th. 3.4]{mastroDB2006} for all $r>0$ the map $A:Z_r(\varphi)\rightarrow Z_r$ is bounded, i.e.
\begin{equation}\label{A-map}
\|Af\|_{Z_r}\le \C\|f\|_{Z_r(\varphi)},\qquad \forall f\in Z_r(\varphi),\quad  r>0, \qquad \C\ne\C(f).
\end{equation}
By means of (\ref{KA}) and (\ref{A-map}), taking into account (\ref{K-inf}) and applying (\ref{f-f'}) with $Kf$ instead of $f$, we easily get
\begin{equation}\label{K-map}
\|Kf\|_{Z_{r+1}}\le \C\|f\|_{Z_r(\varphi)},\qquad \forall f\in Z_r(\varphi), \quad r>0, \qquad \C\ne\C(f),
\end{equation}
i.e. the map $K:Z_r(\varphi)\rightarrow Z_{r+1}$ is bounded for any $r>0$.

Combining this result with the continuous embedding $Z_{r+1}\subset Z_{r+1}(\varphi)$ and the compact embedding $Z_{r+1}(\varphi)\subset Z_{r}(\varphi)$ (cf. $(iii)$ and $(ii)$ respectively), we obtain the following
\begin{theorem}\label{th-K}
For any $r>0$, the map  $K:Z_{r}(\varphi)\rightarrow Z_{r+1}(\varphi)$ is bounded while the map $K: Z_{s}(\varphi)\rightarrow Z_r(\varphi)$ is a compact operator for all $s\ge r$.
\end{theorem}
Finally, we recall that $K$ maps the orthonormal Chebyshev-2th kind polynomials $\{p_n\}_n$ (cf. (\ref{p}) ) according with the following recurrence relation proved in \cite[Corollary 4.3]{CaCrJu97}
\begin{eqnarray}\label{K0}
K p_0 (y)&=&\frac 1 4 \left[\left(2\log 2+\frac 1 2\right)p_0-\frac 1 2p_2(y)\right],\\
\label{K-l}
K p_\ell (y)&=&\frac 1 4 \left[-\frac 1 \ell p_{\ell-2}(y)+\left(\frac 1 \ell+\frac 1 {\ell+2}\right)p_\ell(y)-\frac 1 {\ell+2}p_{\ell+2}(y)\right],\quad \ell >0.
\end{eqnarray}
\subsection{\underline{On the operator $H$}}
The mapping properties of the operator $H$ depend on the smoothness of its bivariate kernel $h(x,y)$. For instance, if $h$ is continuous w.r.t. both the variables, i.e. if $h\in C([-1,1]^2)$ then, similarly to (\ref{K-inf}), we get
\begin{equation}\label{H-inf}
\|Hf\|\le\C\|f\varphi\|,\quad \forall f\in C^0_\varphi,\qquad \C\ne\C(f).
\end{equation}
In view of Theorems \ref{th-D}, \ref{th-K}, we aim to have, for some $s>0$, the map
$H:Z_{s+1}(\varphi)\rightarrow Z_{s}(\varphi)$ is compact or, at least, such that
\begin{equation}\label{hp-H}
\|Hf\|_{Z_{s+1}(\varphi)}\le \C\|f\|_{Z_s(\varphi)},\qquad \forall f\in Z_s(\varphi),\qquad \C\ne\C(f).
\end{equation}
Using the notation $\overline{h}_x$ (resp. $\overline{\overline{h}}_y$) to look at the bivariate function $h(x,y)$ as an univariate function of the second (resp. first) variable, i.e. setting
\begin{equation}\label{hy}
\overline{h}_x(t):=h(x,t) \qquad\mbox{(resp. $\overline{\overline{h}}_y(t):=h(t,y)$)},\qquad -1\le t\le 1,
\end{equation}
a sufficient condition to get (\ref{hp-H}) is given by the following proposition.
\begin{proposition}
\label{prop-H}
\cite[Prop. 4.12]{ JuLu-JCAM97} Let be $s>0$ and $v(x)=(1-x)^\nu(1+x)^\zeta$ with $\nu,\zeta\in [0,1[$. If the kernel $h$ is s.t. $h(x,y)v(x)\varphi(y)\in C^0([-1,1]^2)$ and
$\overline{h}_x v(x)\in Z_s(\varphi)$ uniformly w.r.t. $x\in [-1,1]$, then we have that the map $H:C^0_\varphi\rightarrow Z_s(\varphi)$ is bounded, i.e.
\[
\|Hf\|_{Z_s(\varphi)}\le \C\|f\varphi\|,\quad \forall f\in C^0_\varphi, \qquad \C\ne \C(f).
\]
\end{proposition}
Combining this result with the embedding property $(ii)$, we get, under the hypotheses of Proposition \ref{prop-H}, that
\begin{equation}\label{H-map1}
H:C^0_\varphi\rightarrow Z_r(\varphi) \ \mbox{is a compact operator for all $r\in ]0,s[$}.
\end{equation}
Similarly, recalling $(i)$, under the assumptions of Proposition \ref{prop-H}, we get
\begin{equation}\label{H-map2}
H:Z_r(\varphi)\rightarrow Z_s(\varphi) \ \mbox{is a compact operator for all $r>0$}.
\end{equation}
\subsection{\underline{On the unique solvability of the Prandtl equation}}
Let us write the Prandtl equation (\ref{eq-cont-op}) in the shorter form
\begin{equation}\label{eq-xF}
(D+U)f=g,\qquad U:=\sigma I+K+H
\end{equation}
and suppose that, for a given $s>0$ we have $g\in Z_s(\varphi)$.

In order to prove the existence and uniqueness of the solution $f$ of (\ref{eq-xF}), we use the result in \cite[Cor. 3.8]{Kress}) that follows from the classical Fredholm's alternative theorem.

Taking into account the mapping properties of the operators in (\ref{eq-xF}),
we get
\begin{theorem}\label{uni}
Let us assume that for some $s>0$, we have that $g\in Z_{s}(\varphi)$ and that the map $H: Z_{s+1}(\varphi) \rightarrow Z_s(\varphi)$ is a compact operator. If the map $D+U: Z_{s+1}(\varphi) \rightarrow Z_s(\varphi)$ is injective (i.e. if $Ker(D+U)=\{0\}$ in $Z_{s+1}(\varphi)$) then equation $(D+U)f=g$ admits a unique, stable solution $f^*\in Z_{s+1}(\varphi)$.
\end{theorem}
Recalling Proposition \ref{prop-H}, we also have
\begin{theorem}\label{uni-Marika}
Let us assume that $h(x,y)\varphi(y)\in C^0([-1,1]^2)$ and that $\overline{h}_x\in Z_{s}(\varphi)$ uniformly w.r.t. $x\in [-1,1]$, for some $s>0$. Moreover, let be $0<r\le s$ and suppose that $D+U: Z_{r+1}(\varphi) \rightarrow Z_r(\varphi)$ is an injective map. Then for all $g\in Z_{r}(\varphi)$, the equation $(D+U)f=g$ admits a unique, stable solution $f^*\in Z_{r+1}(\varphi)$.
\end{theorem}
We remark that the previous theorems ensure, besides the unisolvence, also the stability of problem (\ref{eq-xF}), that means the continuity of $(D+U)^{-1}: Z_r(\varphi) \rightarrow Z_{r+1}(\varphi)$, $r\le s$.
\section{Filtered VP interpolation}
For any $n\in\NN$, let us consider the zeros of the $n$--th  Chebyshev--2th kind polynomial $p_n$ (cf. (\ref{p})) and the related Christoffel numbers, respectively given by the formulas
 \[
 x_k:=x_{n,k}=\cos\left( \frac{k\pi}{n+1}\right) \quad\mbox{and}\quad \lambda_k:=\lambda_{n,k}=\frac{\pi}{n+1}\sin^2\left(\frac{k\pi}{n+1}\right),
 \qquad k=1,\ldots, n.
 \]
In alternative to the classical Lagrange interpolation at the nodes $\{x_k\}_k$, we consider the following VP interpolation polynomials \cite{woula_NA}
\begin{equation}\label{VP}
V_n^mf(x):=\sum_{k=1}^n f(x_k)\Phi_{n,k}^m(x), \qquad |x|\le 1,\qquad m<n,
\end{equation}
defined by means of the  so--called fundamental VP polynomials
\begin{equation}\label{fi-filter}
\Phi_{n,k}^m(x)=\lambda_{k}\sum_{j=0}^{n+m-1}\mu_{n,j}^m p_j(x_k)p_j(x),\qquad k=1,\ldots,n,
\end{equation}
where $\mu_{n,j}^m$ are the following VP filtering coefficients
\begin{equation}\label{muj}
\mu_{n,j}^m:=\left\{\begin{array}{ll}
1 & \mbox{if}\quad j=0,\ldots, n-m,\\ [.1in]
\displaystyle\frac{n+m-j}{2m} & \mbox{if}\quad
n-m< j< n+m.
\end{array}\right.
\end{equation}
Equivalently, we can define the fundamental VP polynomials as the following delayed arithmetic means of the Darboux kernels $K_r(x,y):=\sum_{j=0}^rp_j(x)p_j(y)$

\begin{equation}\label{mean-Darboux}
\Phi_{n,k}^m(x)=\frac{\lambda_{k}}{2m}\sum_{r=n-m}^{n+m-1}K_r(x_k, x),\qquad k=1,\ldots,n.
\end{equation}
Moreover,  for  $k=1,\ldots,n$ the following trigonometric form can be found in \cite{filtered1D}
\begin{equation}\label{fi-trig}
\Phi_{n,k}^m(\cos t)=\frac{(-1)^k\sin t_k}{4 m (n+1)}\ \frac{\sin[(n+1)t]}{\sin t}\left[
\frac{\sin[m(t-t_k)]}{\sin^2[(t-t_k)/2]}
-\frac{\sin[m(t+t_k)]}{\sin^2[(t+t_k)/2]}\right],
\end{equation}
where $t_k=\frac{k\pi}{n+1}$.
Similarly to Lagrange polynomials, the interpolation property (see \cite[Section 4]{woula_NA})
\begin{equation}\label{int-VPfund}
\Phi_{n,k}^m(x_h)=\delta_{h,k}=\left\{\begin{array}{ll}
1 & h=k\\
0 & h\ne k
\end{array}\right.\qquad  k=1,\ldots,n,
\end{equation}
holds for any choice of the positive integer $m<n$. Hence we have
\begin{equation}\label{int-VP}
V_n^mf(x_k)=f(x_k),\qquad  k=1,\ldots,n, \qquad m<n.
\end{equation}
Nevertheless, denoted by $\PP_n$ the space of all algebraic polynomials of degree at most $n$, the map $V_n^m: f\rightarrow V_n^mf\in\PP_{n+m-1}$ is not a polynomial projection in the classical sense. It is usually defined as a {\em polynomial quasi-projector}, since its codomain is $\PP_{n+m-1}$ but it preserves all the polynomials of (lower) degree at most $n-m$. Indeed, it has been proved in \cite[Th.4.2]{woula_NA} that $V_n^m$   is a projection onto a non--canonical  polynomial space, the so--called VP space:
\[
\mathrm{S}_n^m:= \textrm{span}\left\{ \Phi_{n,k}^m:\ k=1,\ldots, n\right\},
\]
which is nested between the classical polynomial spaces as follows
\[
\PP_{n-m}\subset S_n^m\subset \PP_{n+m-1}.
\]
We remark that $\dim \mathrm{S}_n^m=\dim \PP_{n-1}=n$. Moreover,  we have \cite[Th.4.2]{woula_NA}
\[
f\in \mathrm{S}_n^m\ \Longleftrightarrow f=V_n^mf.
\]
We recall that, w.r.t. the  scalar product
\[
<f,g>:= \int_{-1}^1f(x)g(x)\varphi(x)dx,
\]
an orthogonal basis of $\mathrm{S}_n^m$ is given by \cite[Th.4.3]{woula_NA}
\begin{equation}\label{q-basis}
\hspace{-.2cm}q_{n,j}^m(x):=\left\{\begin{array}{ll}
p_j(x)&\mbox{if}\quad j=0,\ldots,n\hspace{-.05cm}-\hspace{-.05cm}m,\\ [.1in]
\displaystyle\frac{m\hspace{-.05cm}+\hspace{-.05cm}n\hspace{-.05cm}-
\hspace{-.05cm}j}{2m}p_j(x)-\frac{m\hspace{-.05cm}-
\hspace{-.05cm}n\hspace{-.05cm}+\hspace{-.05cm}j}{2m}p_{2n-j}(x)
&\mbox{if}\quad n\hspace{-.05cm}-\hspace{-.05cm}m<j<n,
\end{array}\right.
\end{equation}
i.e., we have
\begin{equation}\label{span-q}
\mathrm{S}_n^m= \textrm{span} \left\{q_{n,j}^m,\quad j=0,\ldots,n-1 \right\}
\end{equation}
and
\begin{equation}\label{scprod-qj}
<q_{n,i}^m,q_{n,j}^m>=\delta_{i,j} \left\{
\begin{array}{ll}
1 & 0\le j\le n-m,\\
\frac {m^2+(n-j)^2}{2m^2} & n-m<j<n.
\end{array}\right.
\end{equation}
In \cite[Th. 4.3]{woula_NA}  the following bases changes was proved
\begin{equation}\label{nucleo_base_ortogonale}\Phi_{n,k}^m(x)=\lambda_k \sum_{j=0}^{n-1} p_j(x_k)q_{n,j}^m(x).\end{equation}
This formula yields  the following orthogonal expansion
\begin{equation}\label{V_n_ortho}
V_n^m f(x)=\sum_{j=0}^{n-1}c_{n,j}(f)q_{n,j}^m(x),\qquad 0<m<n,
\end{equation}
where
\begin{equation}\label{cj}
c_{n,j}(f)=\sum_{k=1}^{n}\lambda_k p_j(x_k)f(x_k).
\end{equation}
We point out that the coefficients $c_{n,j}(f)$ do not depend on the choice of $m$. Indeed they coincide with the discretization of the Chebyshev--Fourier coefficients
 \[
 c_j(f)=\int_{-1}^1 p_j(y)f(y)\varphi(y), \qquad j=0,\ldots, n-1,
 \]
 by means of the Gauss--Chebyshev quadrature rule based on the zeros of $p_n$.

Taking into account (\ref{inva-D}), in the sequel we will also consider the following modified polynomials
\begin{equation}\label{qtilde-basis}
\hspace{-.2cm}\tilde q_{n,j}^m(x):=\left\{\begin{array}{ll}
\frac{p_j(x)}{j+1}&\mbox{if}\quad j=0,\ldots,n\hspace{-.05cm}-\hspace{-.05cm}m,\\ [.1in]
\displaystyle\frac{m\hspace{-.05cm}+\hspace{-.05cm}n\hspace{-.05cm}-
\hspace{-.05cm}j}{2m}\ \frac{p_j(x)}{j+1}-\frac{m\hspace{-.05cm}-
\hspace{-.05cm}n\hspace{-.05cm}+\hspace{-.05cm}j}{2m}\ \frac{p_{2n-j}(x)}{2n-j+1}
&\mbox{if}\quad n\hspace{-.05cm}-\hspace{-.05cm}m<j<n,
\end{array}\right.
\end{equation}
which yield  an orthogonal basis for the following modified VP space
\[
\tilde S_n^m=\textrm{span} \{\tilde q_{n,j}^m(x):\ j=0,\ldots,n-1\}.
\]
Moreover, we note that
\begin{equation}\label{mix-pr-0}
<q_{n,j}^m, \tilde q_{n,i}^m>=0,\qquad \forall i\ne j
\end{equation}
and, for any $j=1,\ldots,n,$ we have
\begin{equation}\label{mix-pr-1}
<q_{n,j}^m, \tilde q_{n,j}^m>=\left\{\begin{array}{ll}
\frac{1}{j+1}&\mbox{if}\quad j=0,\ldots,n\hspace{-.05cm}-\hspace{-.05cm}m,\\ [.1in]
\displaystyle\frac{(m\hspace{-.05cm}+\hspace{-.05cm}n\hspace{-.05cm}-
\hspace{-.05cm}j)^2}{4m^2(j+1)}+\frac{(m\hspace{-.05cm}-
\hspace{-.05cm}n\hspace{-.05cm}+\hspace{-.05cm}j)^2}{4m^2(2n-j+1)}
&\mbox{if}\quad n\hspace{-.05cm}-\hspace{-.05cm}m<j<n.
\end{array}\right.
\end{equation}
The motivation of considering such modified basis $\tilde q_{n,j}^m$ and space $\tilde S_n^m$, comes from (\ref{inva-D}) which yields
\begin{equation}\label{inva-q}
D\tilde q_{n,j}^m=q_{n,j}^m,\qquad j=0,\ldots,n-1,
\end{equation}
and this identity implies the following
 \begin{proposition}\label{prop}
The map $D: f\in \tilde S_n^m\rightarrow D f\in S_n^m$ is a bijective map and
we have
\[
V_n^m Df=Df,\qquad \forall f\in \tilde S_n^m.
\]
\end{proposition}
About the approximation provided by VP interpolation operators $V_n^m:f\rightarrow V_n^m f$, under the assumption
\begin{equation}\label{ntildem}
m=\theta n, \qquad \mbox{for a fixed $0<\theta<1$} \qquad \quad\mbox{( briefly $n\sim m.$)},
\end{equation}
the following theorem holds (see e.g. \cite{filtered1D, papero-mata}).
\begin{theorem}\label{delavalle}

For any  $n,m\in\NN$ satisfying (\ref{ntildem}), the map
$V_n^m:C^0_\varphi\rightarrow C^0_\varphi$ is bounded, and we have
\begin{equation}\label{delavalle-conv}
\|(f-V_n^mf)\varphi\| \leq \C E_{n-m}(f)_\varphi,\quad \forall f\in C^0_\varphi, \qquad\C\neq \C(n,f).
\end{equation}
Moreover, $\forall r>0$, also the map $V_n^m:Z_r(\varphi)\rightarrow Z_r(\varphi)$ is bounded and the following error estimates holds in H\"older--Zygmund spaces
\begin{equation}\label{delavalle-Zig}
\|f-V_n^mf\|_{Z_r(\varphi)} \leq \frac\C{n^{s-r}}\|f\|_{Z_s(\varphi)},\qquad \forall f\in Z_s(\varphi), \ s\ge r>0, \qquad\C\neq \C(n,f).
\end{equation}
\end{theorem}
Note that for any $n\sim m$, (\ref{delavalle-conv}) and (\ref{Favard-Zig})  yield
\begin{equation}\label{delavalle-inf}
\|(f-V_n^mf)\varphi\| \leq \frac \C{n^r}\|f\|_{Z_r(\varphi)}, \qquad \forall f\in Z_r(\varphi),\ r>0, \qquad  \C\neq \C(n,f).
\end{equation}
In the sequel, we are going to apply the projection map $V_n^m:f\rightarrow V_n^mf\in S_n^m$ in order to get a discrete version of the Prandtl equation (\ref{eq-cont-op}).

Concerning the choice to do on the degree--parameters $n>m$, we remark that $n$ will represent the dimension of the discrete problem while $m$ yields the action ray of VP filter, which we choose according with (\ref{ntildem}). The fixed parameter $0<\theta <1$ acts as a {\it localization parameter} \cite{woula_NA}.

In Section 6 we take, for simplicity, VP filters acting from an even integer $N$ to $2N-1$. In this case, $\theta=\frac 13$,  the action ray is $m=N/2$ and $n=3N/2$ is the number of interpolation nodes that are the Chebyshev zeros of $p_{3N/2}(x)$.

\section{The approximate equation}
Let us assume that the conditions ensuring existence and uniqueness of the solution of the equation  (\ref{iniziale2}) in H\"older--Zygmund spaces are satisfied. In this section we are going to construct a finite dimensional equation to approximate (\ref{iniziale2}). Moreover, we deduce the unique solvability of such approximate equation from that one of original problem (cf. Theorem \ref{stabilita}). Finally we estimate the error we make by approximating the solution of the original problem with that one of the approximated problem.

Let us start by introducing the following discrete approximation of the operator $K$
\begin{equation}\label{def-Kn}
K_n^m f(y):=V_n^m ( K f)(y)=\sum_{j=0}^{n-1}\left[\sum_{k=1}^n\lambda_k p_j(x_k)Kf(x_k)\right]q_{n,j}^m(y).
\end{equation}
By applying Th. \ref{delavalle} and Th. \ref{th-K}, we easily deduce the following result concerning the previous approximation.
\begin{theorem}\label{th-Kn}
For all $r>0$ and any pair of positive integers $n\sim m$, the map $K_n^m:Z_{r}(\varphi)\rightarrow Z_{r+1}(\varphi)$ is bounded and the map $K_n^m:Z_{r}(\varphi)\rightarrow Z_{r}(\varphi)$ is compact. Moreover, for each $f\in Z_{s}(\varphi)$, $s>0$, and any $0<r\le s+1$, we have
\begin{eqnarray}
\label{Kerr-Zig}
\|(K-K_n^m)f\|_{Z_{r}(\varphi)}&\le& \frac \C{n^{s+1-r}}\|f\|_{Z_s(\varphi)},\qquad \C\ne \C(n,m,f),\\
\label{Kerr-inf}
\|(K-K_n^m)f\|_\infty &\le& \frac \C{n^{s+1}}\|f\|_{Z_s(\varphi)},\qquad \C\ne \C(n,m,f).
\end{eqnarray}
\end{theorem}
 In order to approximate the operator $H$, in the definition (\ref{H}), we regard the bivariate kernel function $h$ as the univariate function $\overline{\overline{h}}_y$ (cf. (\ref{hy})), which
we are going to approximate by $V_n^m\overline{\overline{h}}_y$.  Recalling  (\ref{V_n_ortho})--(\ref{cj}), we have
\begin{equation}\label{def-Vhy}
V_n^m\overline{\overline{h}}_y(x)=\sum_{j=0}^{n-1}\left[\sum_{k=1}^n\lambda_kp_j(x_k)h(x_k,y)\right]q_{n,j}^m(x), \qquad |x|\le 1,
\end{equation}
and using in (\ref{H}) this polynomial instead of $\overline{\overline{h}}_y(x)$, we get the operator
\begin{equation}\label{def-Htilde}
\tilde H_n^m f(y):=\frac 1 \pi \int_{-1}^1 V_n^m\overline{\overline{h}}_y(x) f(x)\varphi(x)dx,
\end{equation}
which we approximate as follows
\begin{equation}\label{def-Hn}
H_n^m f(y):=V_n ^m( \tilde H_n^m f)(y)=\sum_{j=0}^{n-1}\left[\sum_{k=1}^n\lambda_k p_j(x_k)\tilde H_n^m f(x_k)\right]q_{n,j}^m(y)  .
\end{equation}
Concerning the previous approximate operators, we prove (see Appendix) the following
\begin{theorem}\label{th-Hn}
Let be $s>0$ and $n\sim m$. If $h(x,y)\varphi(y)\in C^0([-1,1]^2)$ and
$\overline{h}_x\in Z_s(\varphi)$ uniformly w.r.t. $x\in [-1,1]$, then  the maps $\tilde H_n^m:C^0_\varphi\rightarrow Z_s(\varphi)$ and $H_n^m:C^0_\varphi\rightarrow Z_s(\varphi)$ are bounded. Moreover, for each $f\in C^0_\varphi$, and any $0<r\le s$, we have
\begin{eqnarray}
\label{Herr-Zig}
\|(H-H_n^m)f\|_{Z_{r}(\varphi)}&\le& \frac \C{n^{s-r}}\|f\varphi\|,\qquad \C\ne \C(n,m,f),\\
\label{Herr-inf}
\|(H-H_n^m)f\varphi\| &\le& \frac \C{n^{s}}\|f\varphi\|,\qquad \C\ne \C(n,m,f).
\end{eqnarray}
\end{theorem}
By using the previous approximations,  we define the following approximate equation
\begin{equation}\label{eq-approx}
\left(\sigma V_n^m+D+K_n^m+H_n^m\right)\tilde f=V_n^m g.
\end{equation}
For the time being, let us suppose that (\ref{eq-approx}) has a solution $\tilde f$. In this case, we have
\[
D\tilde f=V_n^m g -\sigma V_n^m\tilde f-K_n^m\tilde f-H_n^m\tilde f\in S_n^m,
\]
and recalling Proposition \ref{prop}, the map $D:\tilde S_n^m\rightarrow S_n^m$ is bijective. Hence, we conclude that if a solution of (\ref{eq-approx}) exists, then it belongs to $\tilde S_n^m$.

On the other hand, the fact that $\tilde f\in  \tilde S_n^m$ allows us to look at the approximate equation (\ref{eq-approx}) as the equation resulting by a standard projection method, based on the VP projection $V_n^m :f\rightarrow V_n^mf\in S_m^n$ and applied to the following approximate equation
\begin{equation}\label{eq-project}
\left(\sigma I+D+K+\tilde H_n^m\right)\tilde f= g,
\end{equation}
which is deduced from (\ref{eq-cont-op}) by approximating its last integral.
Of course in the case $h(x,y)\equiv 0$ (see Subsection 5.1) this equation 
 coincides with (\ref{eq-cont-op}) and the same holds if $\overline{\overline{h}}_y\in \PP_{n-m}$ being $\tilde H_n^m= H$ in this case.

Hence, the equation (\ref{eq-approx}), whose solution $\tilde f\in\tilde S_n^m$, can be obtained by applying the VP projection $V_n^m$ at both the sides of the equation (\ref{eq-project}).
Consequently, we can use standard arguments of projection methods in order to deduce  the unique solvability of the approximate equation
\[
T_n^m \tilde f=V_n^m g, \qquad T_n^m:=D+U_n^m,\qquad U_n^m:=\sigma V_n^m +K_n^m+H_n^m,
\]
from the unique solvability of the original Prandtl equation
\[
Tf=g,\qquad T:=D+U,\qquad U:=\sigma I +K+H.
\]
Indeed, by Th. \ref{delavalle}, Th. \ref{th-Kn} and Th. \ref{th-Hn} we get
\begin{equation}\label{U-Un}
\lim_{n\rightarrow+\infty}\|U-U_n^m\|_{Z_r(\varphi)\rightarrow Z_r(\varphi)}=0,\qquad n\sim m,
\end{equation}
where $\|U\|_{A\rightarrow B}$ denotes the norm of the operator map $U:A\rightarrow B$.

Hence, by standard arguments (see e.g. \cite[Th. 4.2]{air}) we can prove the following
\begin{theorem}\label{stabilita}
Let the assumptions of Theorem \ref{uni-Marika} be satisfied. For every $0<r\le s$ and for all $n\sim m$ sufficiently large (say $n>n_0$), the operator $T_n^m:Z_{r+1}(\varphi)\rightarrow Z_{r}(\varphi)$ has a bounded inverse and
\[\sup_{n\sim m} \|(T_n^m)^{-1}\|_{Z_{r}(\varphi)\rightarrow Z_{r+1}(\varphi)}<+\infty.\]
Moreover, the condition number of $T_n^m$ tends to the condition number of $T$,  i.e.
\[
\lim_{\small \begin{array}{c}n\rightarrow \infty\\
n\sim m
\end{array}} \frac{\|T_n^m\|_
{Z_{r+1}(\varphi)\rightarrow Z_{r}(\varphi)}
\|(T_n^m)^{-1}\|_{Z_{r}(\varphi)\rightarrow Z_{r+1}(\varphi)}}{
\|T\|_{Z_{r+1}(\varphi)\rightarrow Z_{r}(\varphi)}\|T^{-1}\|_
{Z_{r}(\varphi)\rightarrow Z_{r+1}(\varphi)}}=1.\]
\end{theorem}
Hence, under the hypotheses of the previous theorem, for any $n\sim m$ sufficiently large and for any $g\in Z_r(\varphi)$, $r>0$, there exists a unique solution $\tilde f_n^m$ of (\ref{eq-approx}), which is a stable polynomial approximation of the solution $f^*\in Z_{r+1}(\varphi)$ of the Prandtl equation (\ref{eq-cont-op}). The following theorem ensures us that as $n\sim m\rightarrow \infty$, we have $\tilde f_n^m\rightarrow f^*$ and the error estimates give the same convergence rate of the error of best approximation of $g$.
\begin{theorem}\label{errore}
Let the assumption of Theorem \ref{uni-Marika} be satisfied. For every $g\in Z_s(\varphi)$, $0<r\le s$ and for all $n\sim m$ sufficiently large (say $n>n_0$), the unique solutions $f^*$ and $\tilde f_n^m$ of (\ref{eq-cont-op}) and (\ref{eq-approx}), respectively, satisfy the following estimates
\begin{eqnarray}
\label{err-Zig}
\|f^*-\tilde f_n^m\|_{Z_r(\varphi)}&\le& \frac{\C}{n^{s-r}}\|g\|_{Z_{s}(\varphi)}, \qquad \C\ne \C(n,m,f^*),\\
\label{err-inf}
\|(f^*-\tilde f_n^m)\varphi\|&\le& \frac{\C}{n^{s}}\|g\|_{Z_{s}(\varphi)}, \qquad \C\ne \C(n,m,f^*).
\end{eqnarray}
\end{theorem}
Comparing this result with the convergence estimates of the methods in  \cite{CaCrJuLu00, DBO}, we succeed in cutting off the typical $\log n$ factors occurring with uniform norms,  reproducing the estimates of the $L^2$ case (see \cite{CaCrJu97}).
\section{On the computation of the approximate solution}
In this section, we suppose that the finite dimensional equation (\ref{eq-approx}) has the unique solution $\tilde f_n^m$ and we face the problem of computing it.

We know  that $\tilde f_n^m\in \tilde S_n^m$, hence choosing a basis in $\tilde S_n^m$, the problem is reduced to the computation of the $n=\dim \tilde S_n^m$ coefficients of  $\tilde f_n^m$ in such basis. This can be done by using the approximate equation (\ref{eq-approx}). Indeed this is an equality in $S_n^m$ and hence, by choosing a basis in $S_n^m$, (\ref{eq-approx}) can be reduced to a system of $n=\dim S_n^m$ equations. In this way, once we have chosen a couple of bases in  $\tilde S_n^m$ and  $S_n^m$, we get a linear system of $n$ equations in $n$ unknowns that are the  coefficients of  $\tilde f_n^m$  in the chosen basis of $\tilde S_n^m$.

Obviously, different choices of the bases give rise to different linear systems to solve. Here we propose to choose the orthogonal bases (\ref{q-basis}) and (\ref{qtilde-basis}) in $S_n^m$ and $\tilde S_n^m$, respectively. We will see that this choice gives us some advantages.

In the sequel,  we adopt simplified notation and  set
\[
q_j:=q_{n,j}^m,\qquad \tilde q_j=\tilde q_{n,j}^m.
\]
Hence we represent the solution $\tilde f:=\tilde f_n^m\in \widetilde{\mathrm{S}}_{n}^m$ of the approximate equation (\ref{eq-approx}) as  follows
\begin{equation}\label{soluzione_apprx}
\tilde f(y)=\sum_{j=0}^{n-1} \tilde f_j \widetilde q_j(y).
\end{equation}
The coefficients vector $ \mathbf{\tilde f}=[\tilde f_0,\ldots, \tilde f_{n-1}]^T$ constitutes the unknown of the linear system we are building.
\begin{lemma}\label{lemmaVn}
We have
\begin{equation}\label{vnf}
\tilde f(y)=\sum_{j=0}^{n-1} \tilde f_j \widetilde q_j(y)\quad \Longrightarrow\quad V_n^m \tilde f(y)= \sum_{j=0}^{n-1}\tilde f_j w_j q_j(y),
\end{equation}
where we set
\begin{equation}\label{wj}
w_j=\left\{\begin{array}{ll}\frac 1 {j+1} & 0\le j \le n-m,\\  [.1in]
\frac 1{2m} \left( \frac{m+n-j}{j+1}+\frac{j-n+m}{2n-j+1}\right) & n-m< j< n .\end{array}\right.
\end{equation}
\end{lemma}
By virtue of this lemma, set
\[
\mathbf{\tilde q}(y):=[\tilde q_0(y),\ldots, \tilde q_{n-1}(y)],\qquad\mathbf{q}(y):=[q_0(y),\ldots, q_{n-1}(y)],\qquad {\mathcal{V}}_n:=\textrm{diag} (w_j)_{j=0,..,n-1},
\]
we have that
\begin{equation}\label{vnf-vec}
\tilde f(y)=\mathbf{\tilde q} (y)\cdot\mathbf{\tilde f}\quad \Longrightarrow\quad V_n^m \tilde f(y)= \mathbf{q} (y)\cdot {\mathcal V}_n\hspace{.1cm} \mathbf{\tilde f}.
\end{equation}
\begin{lemma}\label{lemmaKn}
If $\displaystyle \tilde f(y)=\sum_{j=0}^{n-1} \tilde f_j \widetilde q_j(y)$, then for any $m<n$ we have
\begin{equation}\label{K_n}
K_n^m \tilde f (y)= \sum_{j=0}^{n-1}\left[\alpha_{j+2}\tilde f_{j+2}+\beta_{j}\tilde f_{j}+\gamma_{j-2}\tilde f_{j-2}\right]q_j(y),
\end{equation}
where we set
\begin{eqnarray}
\label{a}
\alpha_\ell &:=&\left\{\begin{array}{ll}
 -\frac 1 {4\ell(\ell+1)}& 2\le \ell\le n-m,\\ [.1in]
-\frac 1 {8m}
{\left[\frac{n+m-\ell}{\ell(\ell+1)}+\frac{\ell-n+m}{(2n-\ell+1)
(2n-\ell+2)}\right]} & n-m<\ell<n,\\
0 & otherwise,
 \end{array} \right.
 \\ [.1in]
 \label{b}
\beta_\ell &:=&\left\{\begin{array}{ll}
\frac 1 4 \left(\frac 12+2\log 2\right) & \ell=0,\\[.07in]
 \frac 1 {2\ell(\ell+2)} & \ell=1,2,\dots,n-m,\\ [.07in]
\frac 1 {4m}\left[\frac{n+m-\ell}{\ell(\ell+2)}+\frac{\ell-n+m}{(2n-\ell)
(2n-\ell+2)}\right]& n-m<\ell\le n-2,\\ [.07in]
\frac 1 {8 m}\left[
\frac{(m+1)(3n-1)}{n(n^2-1)}+\frac{(m-1)(3n+7)}{(n+1)(n+2)(n+3)}\right]& \ell= n-1,\end{array}\right.
\\ [.1in]
\label{c}
\gamma_\ell &:=&\begin{cases}
-\frac 1 8 & \ell=0,\\  -\frac 1 {4(\ell+1)(\ell+2)} & \ell=1,2,\dots,n-m,\\
-\frac 1 {8m}\left[\frac{n+m-\ell}{(\ell+1)(\ell+2)}+\frac{\ell-n+m}{(2n-\ell+1)
(2n-\ell)}\right]& \ell=n-m+1,\dots,n-3,\\
0 &  otherwise. \end{cases}
\end{eqnarray}
\end{lemma}
In vector form, Lemma \ref{lemmaKn} states that
\begin{equation}\label{Knf-vec}
\tilde f(y)=\mathbf{\tilde q} (y)\cdot\mathbf{\tilde f}\quad \Longrightarrow\quad K_n^m \tilde f(y)= \mathbf{q} (y)\cdot {\mathcal A}_n\hspace{.1cm} \mathbf{\tilde f},
\end{equation}
where we set
{\small{\begin{equation}\label{Mat-An}
\mathcal{A}_n :=
\begin{pmatrix}
\beta_0 & 0 & \alpha_2 &  &  &   &  & \\
0 & \beta_1 & 0 & \alpha_3 &  &  & \textrm{\huge 0} & \\
\gamma_0 & 0 &  \beta_2 & 0 & \alpha_4 & &  &  \\
 & \gamma_1 & 0 &  \beta_3 & 0 & \alpha_5  &  &  \\
 &  & \ddots & \ddots & \ddots & \ddots &  \ddots & \\
 &  &  & \ddots & \ddots & \ddots & \ddots & \alpha_{n-1}\\
 &  & \textrm{\huge 0} &  & \ddots & \ddots & \ddots &  0\\
 &  &  &  &  & \gamma_{n-3} & 0 & \beta_{n-1}
\end{pmatrix}.
\end{equation}}}
\begin{lemma}\label{lemmaH_n}
Set $Q_i:=<q_i,\widetilde q_i>_\varphi$, $i=0,\ldots, n-1$ (cf. (\ref{mix-pr-1})). If $\tilde f$ is given by (\ref{soluzione_apprx}) then
\begin{equation}\label{sistema_Hn}
 H_n^m \tilde f (y)=\sum_{j=0}^{n-1}\left[\frac 1 \pi\sum_{i=0}^{n-1}\sum_{s=1}^{n}\sum_{k=1}^{n}\lambda_k \lambda_s p_i(x_s)p_j(x_k)h(x_s,x_k)\tilde f_i Q_i\right] q_j(y).
 \end{equation}
\end{lemma}
In vector form, Lemma \ref{lemmaH_n} states that
\begin{equation}\label{Hnf-vec}
\tilde f(y)=\mathbf{\tilde q} (y)\cdot\mathbf{\tilde f}\quad \Longrightarrow\quad H_n^m \tilde f(y)
= \displaystyle \mathbf{q} (y)\cdot {\mathcal B}_n\hspace{.1cm} \mathbf{\tilde f},
\end{equation}
with
\[
{\mathcal B}_n:= \frac 1\pi \left({\mathcal{P}}_n\Lambda_n\right){\mathcal H}_n\left({\mathcal{P}}_n\Lambda_n\right)^T {\mathcal Q}_n,
\]
where $\Lambda_n$ and $ {\mathcal{Q}_n}$ are the following diagonal matrices
\[
\Lambda_n:=\textrm{diag} (\lambda_j)_{j=0,..,n-1},\qquad {\mathcal{Q}_n}:=\textrm{diag} (<q_j,\tilde q_j>)_{j=0,..,n-1},
\]
and $\mathcal{P}_n, \mathcal{H}_n$ are the following square matrices of order $n$
$$\mathcal{P}_n(i,j)=p_{i-1}(x_j),
\qquad \mathcal{H}_n(i,j)=h(x_j,x_i),\qquad 1\le i, j\le n. $$
Now let us write the approximate equation (\ref{eq-approx}) as follows
\begin{equation}\label{eq-vec}
(D+U_n^m)\tilde f(y)=V_n^mg(y),\qquad U_n^m:=\sigma V_n^m+K_n^m+H_n^m.
\end{equation}
The previous lemmas and (\ref{inva-q}) ensure that
\[
\tilde f(y)=\mathbf{\tilde q} (y)\cdot\mathbf{\tilde f}\quad \Longrightarrow\quad
\left\{\begin{array}{rl}
\displaystyle U_n^m \tilde f(y)
&= \displaystyle \mathbf{q} (y)\cdot \left(\sigma {\mathcal V}_n+{\mathcal A}_n+{\mathcal B}_n\right)\hspace{.1cm} \mathbf{\tilde f},
\\[.1in] \displaystyle
D \tilde f(y)
&= \displaystyle \mathbf{q} (y)\cdot  \mathbf{\tilde f}.
\end{array}\right.
\]
Hence, taking into account that
\begin{equation}\label{g-vec}
V_n^mg(y)=\mathbf{q} (y)\cdot  \mathbf{g},\qquad
\mathbf{g}:=\left[g_0,\ldots, g_{n-1}\right]^T,
\qquad g_j:=c_{n,j}(g)=\sum_{k=1}^n\lambda_k p_j(x_k) g(x_k),
\end{equation}
we rewrite (\ref{eq-vec}) as follows
\[
\mathbf{q} (y)\cdot \left({\mathcal I}_n+\sigma {\mathcal V}_n+{\mathcal A}_n+{\mathcal B}_n\right)\hspace{.1cm} \mathbf{\tilde f}= \mathbf{q} (y)\cdot \mathbf{g}.
\]
Consequently, the unknown vector $\mathbf{\tilde f}$ results to be the unique solution of the following linear system
\begin{equation}\label{sist}
\left({\mathcal I}_n+\sigma {\mathcal V}_n+{\mathcal A}_n+{\mathcal B}_n\right)\hspace{.1cm} \mathbf{\tilde f}= \mathbf{g}
\end{equation}
which yields the unique solution
$\tilde f(y)=\mathbf{\tilde q} (y)\cdot\mathbf{\tilde f}$ of the approximate equation (\ref{eq-vec}).
\subsection{The special case $H\equiv 0$}
The case $H\equiv 0$ is a particularly favorite case since, in (\ref{sist}), ${\mathcal B_n}$ is the null matrix and the linear system (\ref{sist}) becomes
\begin{equation}\label{sist-H0}
{\mathcal A_n^\prime}\hspace{.1cm} \mathbf{\tilde f}= \mathbf{g}, \qquad
{\mathcal A_n^\prime}:= {\mathcal I}_n+\sigma {\mathcal V}_n+{\mathcal A}_n.
\end{equation}
The matrix ${\mathcal A_n^\prime}$ of such system has not null only the $(i,j)$--entries s.t. $|i-j|\in\{0,2\}$. It differs from the matrix ${\mathcal A_n}$ in (\ref{Mat-An}) only for the diagonal entries that are given by
\begin{equation}\label{delta-j}
\delta_{j,j}:={\mathcal A_n^\prime}(j,j)=1+\sigma w_j+\beta_j\qquad j=0,\ldots, n-1.
\end{equation}
The matrix  ${\mathcal A_n^\prime}$ is strictly diagonally dominant, since by (\ref{wj}) and (\ref{a})--(\ref{c}) it easily follows that for any $\sigma\in \RR$, we have
$$|1+\sigma w_0+\beta_0|> |\alpha_{2}|, \quad |1+\sigma w_1+\beta_1|> |\alpha_{3}|,$$
$$|1+\sigma w_j+\beta_j|>|\gamma_{j-2}|+|\alpha_{j+2}|, \quad 2\le j\le n-3,$$
$$|1+\sigma w_{n-2}+\beta_{n-2}|> |\gamma_{n-3}|, \quad |1+\sigma w_{n-1}+\beta_{n-1}|> |\gamma_{n-2}|.$$
As a consequence, the Gaussian elimination method for solving the system ${\mathcal A_n^\prime} \widetilde{\mathbf{f}}_n=\mathbf{g}$ can be applied avoiding pivoting strategy to ensure the stability of the algorithm (see e.g.\cite[\S 3.4.10]{golub}). Therefore,   the 2-bandwidth of the matrix is preserved and the Gaussian elimination requires $5n$  long operations.

The condition number of the matrix ${\mathcal A^\prime_n}$ is defined in the usual way
 \[
\kappa({\mathcal A}^\prime_n):=\|\mathcal{A}^\prime_n\|\|\left(\mathcal{A'}_n\right)^{-1}\|,
 \]
 where $\|\cdot\|$ denotes any matrix norm satisfying the submultiplicative property
 (i.e. $\|A B\|\le \|A\| \|B\|$).
 The following theorem proves that $\kappa({\mathcal A}^\prime_n)$ tends to a finite limit as the dimension of the system $n\rightarrow +\infty.$
\begin{theorem}\label{th-kappa}
Under the previous setting, there exists  the limit
$$\lim_{n\rightarrow + \infty} \kappa(\mathcal{A'}_n)<\infty .$$
\end{theorem}
We recall that the  linear system conditioning has been studied also for the method in \cite{DBO} where the authors, in a more general context, proved that the condition numbers diverge as $n\log^3 n$.
\section{Numerical tests}
For simplicity, in order to test the behavior of the numerical method for increasing values of $n\sim m$, we take increasing values of $N\in\NN$ even and focus on the particular case
\begin{equation}\label{N}
m=\frac N 2,\qquad\mbox{and}\qquad n=\frac 32 N,\qquad \mbox{with $N\in\NN$ even}.
\end{equation}
This corresponds to take $\theta=1/3$ in $m=\theta n$, and VP means from $N=(n-m)$ to $2N-1=(n+m)-1$.

We are going to show the performance of the proposed method by some numerical experiments, making comparisons with the method introduced in \cite{DBO} and the collocation--quadrature method studied in \cite{CaCrJu97, CaCrJuLu00},  all of them based on the Lagrange interpolation at the same nodes. Moreover, we also make a comparison  with  the method proposed in \cite{calio} (cf. Example 6.4).

In all the cases the better performance of our method is displayed.

Although the well-conditioning of our linear systems has been proved only in the special case $H\equiv 0$, the numerical evidence shows the uniform boundedness of the condition numbers also in the general case.

For each example, denoted by $X$ a sufficiently large uniform mesh of $[-1,1]$,  we compute the absolute weighted errors
\begin{equation}\label{Err-test}
\mathcal{E}_n:=\max_{x\in X}\biggm(|f(x)-f_n(x)|\varphi(x)\biggm),
\end{equation}
where $n$ indicates the number of the collocation nodes and $f_n$ denotes the numerical solution of the applied method. Moreover, we compute the condition numbers $\cond_n$, w.r.t. the infinity matrix norm, of the involved linear systems of dimension $n$.

In order to distinguish the implemented numerical methods, we use the superscripts $VP$, \cite{CaCrJuLu00}, \cite{DBO} and \cite{calio} in order to indicate our method and the ones in \cite{CaCrJu97, CaCrJuLu00}, \cite{DBO}, \cite{calio}, respectively.

We point out that all the computations were performed in $16-$digits arithmetic and  the solutions of the linear systems have been computed by Gaussian elimination method.

\begin{es}\label{es1} Consider the integral equation
\[
(I+D+H)f=g,
\]
where
\begin{eqnarray*}
 h(x,y)&=&x(y^2|y|+x|x|),\\
g(y)&=&y\left[\left(1+\frac{4y}{15\pi}\right)|y|+\frac{6}{\pi}+\frac{3y^2-2}{\pi\sqrt{1-y^2}}\log\left(\frac{1+\sqrt{1-y^2}}{1-\sqrt{1-y^2}}
\right)\right].
\end{eqnarray*}
This example  can be  found in \cite{CaCrJu97} and the exact solution is $f(y)=y|y|$. In Table \ref{tab:b1} we display the results obtained implementing our method and those in \cite{DBO}, \cite{CaCrJu97, CaCrJuLu00}.

\begin{table}[h]
\caption{Example \ref{es1}}
\label{tab:b1}
\begin{center}
\begin{tabular}{|c|l|l|l|l|l|l|} \hline
$n$  & $\cond_n^{VP}$ & $\mathcal{E}_n^{VP}$ & $\cond_n^{\mbox{\tiny\cite{DBO}}}$ & $\mathcal{E}_n^{\mbox{\tiny\cite{DBO}}}$ & $\cond_n^{\mbox{\tiny\cite{CaCrJuLu00}}}$ & $\mathcal{E}_n^{\mbox{\tiny\cite{CaCrJuLu00}}}$\\ \hline
$8$  & $1.95$&  $2.91e-3 $& $4.89$     & $8.50e-3 $ & $4.36$     & $3.13 $  \\ \hline
$16$ & $2.05$&  $8.08e-4 $& $9.12 $    & $2.27e-3 $ & $8.24 $    & $8.26e-3 $\\ \hline
$32$ & $2.10$&  $2.13e-4 $& $1.75e+1$  & $6.03e-4 $ & $1.59e+1 $    & $2.15e-3 $\\ \hline
$64$ & $2.12$&  $5.48e-4 $& $3.44e+1$  & $1.56e-4 $  & $3.14e+1 $    & $5.52e-4 $\\ \hline
$128$& $2.13$&  $1.34e-5 $& $6.81e+1$  & $3.94e-5 $ & $6.23e+1 $    & $1.40e-4 $\\ \hline
$256$& $2.14$&  $2.83e-6 $& $1.35e+2$  & $9.24e-6 $& $1.24e+2 $    & $3.53e-5 $\\ \hline
$512$& $2.15$&  $3.81e-7 $& $2.70e+2$  & $2.54e-6 $& $2.47e+2 $    & $8.86e-6 $\\ \hline
\end{tabular}
\end{center}
\end{table}
\end{es}
\begin{es}\label{es2} Consider the integral equation (\ref{iniziale2}) with
$$\sigma=1,\quad h(x,y)=\left|\cos\left(y-\frac \pi 4\right)\right|^\frac 92+|\sin(x)|^\frac 72, \quad \ g(y)=|y|^{\frac{11}2}.$$
Since in this case the exact solution is unknown, the errors showed in Table \ref{tab:b2}  have been computed according to (\ref{Err-test}) but replacing $f$ by $f_{n}^{VP}$ with $n=1024$.

\begin{table}[h]
\caption{Example \ref{es2}}
\label{tab:b2}
\begin{center}
\begin{tabular}{|c|l|l|l|l|l|l|} \hline
$n$  & $\cond_n^{VP}$ & $\mathcal{E}_n^{VP}$ & $\cond_n^{\mbox{\tiny\cite{DBO}}}$& $\mathcal{E}_n^{\mbox{\tiny\cite{DBO}}}$ & $\cond_n^{\mbox{\tiny\cite{CaCrJuLu00}}}$& $\mathcal{E}_n^{\mbox{\tiny\cite{CaCrJuLu00}}}$\\ \hline
$8$  & $2.47$&  $6.91e-6  $& $6.13$     & $9.19e-5 $  & $5.48$     & $3.78e-04$  \\ \hline
$16$ & $2.60$&  $8.50e-8  $& $1.16e+1 $ & $9.22e-7 $  & $1.02e+1$  & $1.44e-06$   \\ \hline
$32$ & $2.66$&  $1.08e-9  $& $2.24e+1$  & $2.18e-8 $  & $1.98e+1$  & $2.18e-08$   \\ \hline
$64$ & $2.69$&  $1.39e-11 $& $4.39e+1$  & $7.97e-10$  & $3.91e+1$  & $7.97e-10$ \\ \hline
$128$& $2.71$&  $2.03e-13 $& $8.67e+1$  & $3.48e-11$  & $7.74e+1$  & $3.48e-11$ \\ \hline
$256$& $2.72$&  $4.53e-15 $& $1.72e+2$  & $1.48e-12$  & $1.54e+2$  & $1.63e-11$  \\ \hline
\end{tabular}
\end{center}
\end{table}
\end{es}

\begin{es}\label{es3} Consider the integral equation (\ref{iniziale2}) with
$$\sigma=1,\quad h(x,y)=|y|+|x|,$$
$$g(y)=2+\frac{|y|}2+\frac 2 {3\pi}+\frac 1 4 (1-2y^2+\log 4),
$$
whose exact solution is $f(y)=1.$  We compare our method with those  in \cite{DBO} and \cite{CaCrJu97, CaCrJuLu00}, reporting in the Table \ref{tab:b3} the results.

\begin{table}[h]
\caption{Example \ref{es3}}
\label{tab:b3}
\begin{center}
\begin{tabular}{|c|l|l|l|l|l|l|} \hline
$n$  & $\cond^{VP}$ & $\mathcal{E}_n^{VP}$ & $\cond^{\mbox{\tiny\cite{DBO}}}$ & $\mathcal{E}_n^{\mbox{\tiny\cite{DBO}}}$ & $\cond^{\mbox{\tiny\cite{CaCrJuLu00}}}$ & $\mathcal{E}_n^{\mbox{\tiny\cite{CaCrJuLu00}}}$\\ \hline
$8$  & $2.64$&  $5.44e-4 $& $5.48$     & $1.82e-3$ & $5.01$     & $1.82e-3$  \\ \hline
$16$ & $2.77$&  $1.46e-4 $& $1.05e+1$    & $5.05e-4 $ & $9.50$     & $5.05e-4$  \\ \hline
$32$ & $2.84$&  $3.80e-5 $& $2.04e+1$  & $1.33e-4 $  & $1.84e+1$     & $1.33e-4$ \\ \hline
$64$ & $2.87$&  $9.70e-6 $& $4.00e+1$  & $ 3.44e-5$ & $3.62e+1$     & $ 3.44e-5$   \\ \hline
$128$& $2.89$&  $2.45e-6 $& $7.91e+1$  & $8.74e-6$ & $7.19e+1$     & $8.74e-6$ \\ \hline
$256$& $2.89$&  $6.15e-7 $& $1.57e+2$  & $ 2.20e-6$ & $1.43e+2$     & $ 2.20e-6$ \\ \hline
$512$& $2.90$&  $1.54e-7 $& $3.13e+2$  & $5.53e-7$ & $2.85e+2$     & $5.53e-7$ \\ \hline
\end{tabular}
\end{center}
\end{table}
\end{es}

\begin{es} \label{es4} Consider the integral equation 
\[
(2I+D+K)f=g,
\]
with  $g(x)$ such that the exact solution is  $f(x)=\sqrt{(1-x^2)^3}$.
This example can be found in \cite{calio} where the proposed method provides an approximation of $f$ with at most $2$ exact decimal digits. In Table \ref{tab:b4}, besides the results of our method also the results of the methods in \cite{DBO}, \cite{CaCrJu97, CaCrJuLu00} are presented.

In this case $H\equiv 0$ and, as pointed up in Subsection 5.1, our method leads to solve the system (\ref{sist-H0}) in about $5n$ long operations, whereas the linear systems associated with the methods in \cite{DBO}, \cite{CaCrJu97, CaCrJuLu00} have a full matrix and can be solved by Gaussian elimination method with pivoting, requiring at least $n^3/3$ long operation.
\begin{table}[h]
\caption{Example 4}
\label{tab:b4}
\centering
\begin{center}
\begin{tabular}{|c|l|l|l|l|l|l|} \hline
$n$ & $\cond^{VP}$ & $\mathcal{E}_n^{VP}$ & $\cond^{\mbox{\tiny\cite{DBO}}}$  & $\mathcal{E}_n^{\mbox{\tiny\cite{DBO}}}$ & $\cond^{\mbox{\tiny\cite{CaCrJuLu00}}}$ & $\mathcal{E}_n^{\mbox{\tiny\cite{CaCrJuLu00}}}$\\ \hline
$8 $   & $2.93$ &  $3.24e-4 $ & $4.99    $ &$1.53e-3 $ & $5.24    $&$1.57e-2$ \\ \hline
$16$   & $3.22$ & $1.97e-5 $  & $9.01    $ &$1.01e-4 $ & $1.00e+1$ &$1.97e-3$ \\ \hline
$32$   & $3.36$ & $5.87e-7 $  & $1.68e+1 $ &$2.56e-6 $ & $1.93e+1$ &$2.59e-4$ \\ \hline
$64$   & $3.43$&  $1.84e-8 $  & $3.24e+1 $ &$9.67e-8 $ & $3.77e+1$ &$3.35e-5$  \\ \hline
$128$  & $3.47$&  $6.05e-10 $ & $6.35e+1 $ &$1.57e-9 $ & $7.42e+1$ &$4.28e-6$ \\ \hline
$256$  & $3.49$&  $4.69e-12 $ & $1.25e+2 $ &$4.94e-11 $& $1.47e+2$ &$5.41e-7 $ \\ \hline
$512$  & $3.50$&  $7.41e-14 $ & $2.50e+2 $ &$1.99e-12 $& $2.93e+2$ &$6.80e-8 $ \\ \hline
\end{tabular}
\end{center}
\end{table}
\end{es}

\section{Appendix}
\subsection{Proof of Theorem \ref{th-Hn}}

Firstly note that, by (\ref{def-Htilde}) and (\ref{VP}), we have
\[
\tilde H_n^mf(y)=\int_{-1}^1 {\cal H} (x,y)f(x)\varphi(x)dx, \qquad {\cal H}(x,y):=\sum_{k=1}^n h(x_k,y)\Phi_{n,k}^m(x).
\]
Let us prove that the boundedness of $\tilde H_n^m:C^0_\varphi\rightarrow Z_s(\varphi)$ follows from Proposition \ref{prop-H}. Indeed, we recall that (cf. Theorem \ref{delavalle})
\begin{equation}\label{leb}
\|V_n^m\|_{C^0_\varphi\rightarrow C^0_\varphi}:=\sup_{\|f\varphi\|\le 1}\|(V_n^mf)\varphi\|=\sup_{|x|\le 1}\sum_{k=1}^n\frac{|\Phi_{n,k}^m(x)|\varphi(x)}{\varphi(x_k)}\le \C\ne\C(n,m).
\end{equation}
Hence, from the assumptions on $h$, we deduce that
\begin{itemize}
\item [(a)] ${\cal H}(x,y)\varphi(x)\varphi(y)$ is a continuous function w.r.t. both the variables $x,y\in [-1,1]$
\item [(b)] For all $n\in\NN$ and $|x|\le 1$, taking into account that
\[
E_n(\overline{h}_{x_k})_\varphi= \|(\overline{h}_{x_k}-P^*_{x_k})\varphi\|\le \frac \C{n^s}, \quad k=1,\ldots, n, \qquad \C\ne\C (n,k),
\]
we set
\[
P^*(y):=\sum_{k=1}^n P^*_{x_k}(y)\Phi_{n,k}^m(x)\varphi(x),\qquad |y|\le 1.
\]
Hence we have
\begin{eqnarray*}
E_n(\overline{{\cal H}}_x\varphi(x))_\varphi&=&\inf_{P\in\PP_n} \sup_{|y|\le 1}
\left| {\cal H}(x,y)\varphi(x)-P(y)\right|\varphi(y)\\
&\le&\sup_{|y|\le 1}
\left| {\cal H}(x,y)\varphi(x)-P^*(y)\right|\varphi(y)\\
&=& \sup_{|y|\le 1}\left|\sum_{k=1}^n \left[ h(x_k,y)-P^*_{x_k}(y)\right]\Phi_{n,k}^m(x)\varphi(x)\right|\varphi(y)\\
&\le& \sum_{k=1}^n\frac{|\Phi_{n,k}^m(x)|\varphi(x)}{\varphi(x_k)}
\sup_{|y|\le 1}\left|\overline{h}_{x_k}(y)-P^*_{x_k}(y)\right|\varphi(y)\\
&\le& \frac C{n^s}\sum_{k=1}^n\frac{|\Phi_{n,k}^m(x)|\varphi(x)}{\varphi(x_k)}
\le \frac \C{n^s}, \qquad \C\ne \C(n,x).
\end{eqnarray*}
\end{itemize}
In conclusion, by virtue of (a) and (b), the kernel ${\cal H}$ satisfies the assumptions of Proposition \ref{prop-H} with $v=\varphi$,  so that the operator ddefined by this kernel, namely $\tilde H_n^m:C^0_\varphi\rightarrow Z_s(\varphi)$ is bounded.

Consequently, the map $H_n^m:C^0_\varphi\rightarrow Z_s(\varphi)$ is bounded too, since it is the following composition of bounded maps (cf. Theorem \ref{delavalle})
\[
 H_n^m:C^0_\varphi\begin{array}{c} \mbox{\scriptsize $\tilde H_n^m$} \\ [-.1in] \longrightarrow\end{array} Z_s(\varphi)\begin{array}{c} \mbox{\scriptsize $V_n^m$} \\ [-.1in] \longrightarrow\end{array} Z_s(\varphi).
\]
In order to prove (\ref{Herr-inf}), let us arbitrarily fix $f\in C^0_\varphi$ and $|y|\le 1$.\newline By using (\ref{VP}), we note that
\[
H_n^mf(y)=\sum_{j=1}^n\tilde H_n^mf(x_j)\Phi_{n,j}^m(y)=
\int_{-1}^1\left(\sum_{j=1}^n\sum_{k=1}^n h(x_k,x_j)\Phi_{n,k}^m(x)\Phi_{n,j}^m(y)\right)f(x)\varphi(x)dx,
\]
that means $H_n^mf$ can be obtained by replacing $h$ with its bivariate VP interpolation polynomial based on tensor--product Chebyshev nodes of the second kind (see \cite{lncs_nostro, paperoAMC}).

Hence, we write
\begin{eqnarray*}
Hf(y)-H_n^mf(y)&=&\int_{-1}^1\left(h(x,y)-\sum_{j=1}^nh(x,x_j)\Phi_{n,j}^m(y)\right)f(x)\varphi(x) dx\\
&+&\int_{-1}^1\sum_{j=1}^n\Phi_{n,j}^m(y)\left(h(x,x_j)-
\sum_{k=1}^nh(x_k,x_j)\Phi_{n,k}^m(x)\right)f(x)\varphi(x)dx.
\end{eqnarray*}
Then,  using (\ref{VP}) (\ref{hy}) and (cf. (\ref{int-VP}))
\begin{equation}\label{sum-fund}
\sum_{k=1}^n\Phi_{n,k}^m(x)=1,\qquad \forall |x|\le 1,
\end{equation}
we get
\begin{eqnarray*}
Hf(y)-H_n^mf(y)&=& \int_{-1}^1\bigg(\overline{h}_x(y)-V_n^m (\overline{h}_x)(y)\bigg)f(x)\varphi(x) dx\\
&+&\int_{-1}^1\sum_{j=1}^n\Phi_{n,j}^m(y)\sum_{k=1}^n\Phi_{n,k}^m(x)
\bigg(\overline{h}_x(x_j)-\overline{h}_{x_k}(x_j)\bigg)f(x)\varphi(x)dx.
\end{eqnarray*}
Consequently
\begin{eqnarray*}
&&|Hf(y)-H_n^mf(y)|\varphi(y)\le\int_{-1}^1\bigg|\overline{h}_x(y)-V_n^m (\overline{h}_x)(y)\bigg| \varphi(y)|f(x)|\varphi(x) dx\\
&+&\sum_{j=1}^n |\Phi_{n,j}^m(y)|\varphi(y)\int_{-1}^1\sum_{k=1}^n|\Phi_{n,k}^m(x)|
\bigg|\overline{h}_x(x_j)-\overline{h}_{x_k}(x_j)\bigg| |f(x)|\varphi(x)dx\\
&=:&A+B.
\end{eqnarray*}
On the other hand, the hypothesis  $\overline{h}_x\in Z_s(\varphi)$ uniformly w.r.t. $x\in [-1,1]$ implies that
$\forall n\in\NN$ and $\forall |x|\le 1$ we have
\begin{equation}\label{eq1}
E_n(\overline{h}_x)_\varphi\le \frac \C{n^s},\qquad  \C\ne \C(n,x),
\end{equation}
and in particular, there exists $P^*\in\PP_n$ such that
\begin{equation}\label{eq2}
\sup_{|y|\le 1}\bigg|\overline{h}_x(y)-P^*(y)\bigg|\varphi(y)\le \frac \C{n^s},\qquad  \C\ne \C(n,x).
\end{equation}
Hence, concerning $A$, by (\ref{delavalle-conv}) and (\ref{eq1}), recalling that $m=\theta n$ with a fixed $0<\theta<1$, we get
\[
A\le \int_{-1}^1 E_{n-m}(\overline{h}_x)_\varphi |f(x)|\varphi(x) dx\le \frac\C{(n-m)^s} \int_{-1}^1 |f(x)|\varphi(x) dx\le \frac\C{n^s}\|f\varphi\|.
\]
Moreover, regarding $B$, by (\ref{eq2}) and (\ref{leb}), we have
\begin{eqnarray*}
B\hspace{-.2cm}&\le&\hspace{-.2cm}\|f\varphi\| \sum_{j=1}^n\frac{|\Phi_{n,j}^m(y)|\varphi(y)}{\varphi(x_j)}
\int_{-1}^1\sum_{k=1}^n\frac{|\Phi_{n,k}^m(x)|\varphi(x)}{\varphi(x_k)}
\bigg|(\overline{h}_x\hspace{-.1cm}-\hspace{-.1cm}P^*)(x_j)-
(\overline{h}_{x_k}\hspace{-.1cm}-\hspace{-.1cm}P^*)(x_j)\bigg| \frac{\varphi(x_j)}{\varphi(x)}dx\\
&\le&\hspace{-.2cm}\frac \C{n^s} \|f\varphi\| \sum_{j=1}^n\frac{|\Phi_{n,j}^m(y)|\varphi(y)}{\varphi(x_j)}
\int_{-1}^1\sum_{k=1}^n\frac{|\Phi_{n,k}^m(x)|\varphi(x)}{\varphi(x_k)}
\frac{dx}{\varphi(x)}\\
&\le&\hspace{-.2cm} \frac \C{n^s} \|f\varphi\|\int_{-1}^1
\frac{dx}{\varphi(x)}\le\frac \C{n^s}\|f\varphi\|, \qquad\mbox{with $\C\ne\C(n,f,y)$}.
\end{eqnarray*}
Finally, let us prove (\ref{Herr-Zig}). Taking into account that $H_n^mf\in\PP_{n+m-1}$ (cf. Proposition \ref{prop-H}) and recalling that $H:C^0_\varphi\rightarrow Z_s(\varphi)$ is bounded, for any $f\in C^0_\varphi$ we note that
\begin{eqnarray*}\sup_{k\ge n+m-1} (k+1)^r E_k(Hf-H_n^mf)_{\varphi}&=&\sup_{k\ge n+m-1} (k+1)^r E_k(Hf)_{\varphi}\\
&\le& \C \sup_{k\ge n+m-1} \frac{\|Hf\|_{Z_s(\varphi)}}{k^{s-r}}\le \C\frac{\|f\varphi\|}{n^{s-r}}.
\end{eqnarray*}
Moreover, by (\ref{Herr-inf}) we get
\begin{eqnarray*}
\sup_{k< n+m-1} (k+1)^r E_k(Hf-H_n^mf)_{\varphi}&\le& \|(Hf-H_n^mf)\varphi\|\sup_{k< n+m-1} (k+1)^r \\
&\le& \C \frac{\|f\varphi\|}{n^s}(n+m)^r\le \C\frac{\|f\varphi\|}{n^{s-r}}.
\end{eqnarray*}
Consequently, by (\ref{Herr-inf}) we have
\[
\|(H-H_n^m)f \|_{Z_{r}(\varphi)}=\| (Hf-H_n^mf)\varphi \| + \sup_{k\in\NN} (k+1)^{r} E_k(Hf-H_n^mf)_{\varphi}\le \C\frac{\|f\varphi\|}{n^{s-r}}
\]
and (\ref{Herr-Zig}) follows.
\subsection{Proof of Theorem \ref{errore}.}
Note that we can write
\begin{eqnarray*}
f^*-\tilde f_n^m&=&(D+U_n^m)^{-1}\left[(D+U_n^m)f^*- V_n^mg \right]\\
&=&(D+U_n^m)^{-1}\left[(g- V_n^mg)+ (U_n^m-U)f^*\right].
\end{eqnarray*}
Hence, by the uniform boudedness of $(D+U_n^m)^{-1}:Z_r(\varphi)\rightarrow Z_{r+1}(\varphi)$ (cf. Th. \ref{stabilita}), for any $0< r\le s$ we deduce
\[
\|f^*-\tilde f_n^m\|_{Z_{r+1}(\varphi)}\le \C \left[\|g- V_n^mg\|_{Z_r(\varphi)}+ \|(U_n^m-U)f^*\|_{Z_{r}(\varphi)}\right]
\]
and the statement follows from Th. \ref{delavalle} and (\ref{U-Un}), taking into account that, by hypothesis, $g\in Z_s(\varphi)$, $f^*\in Z_{s+1}(\varphi)$ and  $\|g\|_{Z_s(\varphi)}\sim \|f^*\|_{Z_{s+1}(\varphi)}$.
\subsection{Proof of Lemma \ref{lemmaVn}.}
Let be $\displaystyle\tilde f(y)=\sum_{j=0}^{n-1}\tilde f_j \tilde q_j(y)$. Recalling (\ref{V_n_ortho})--(\ref{cj}), we have
$\displaystyle V_n^m\tilde  f(y)= \sum_{r=0}^{n-1}q_r(y)c_{n,r}(\tilde f),$
where
\begin{eqnarray*}
c_{n,r}(\tilde f)&=&\sum_{k=1}^{n}\lambda_k p_r(x_k)\tilde f(x_k)
\\
&=&\left(\sum_{j=0}^{n-m}+ \sum_{j=n-m+1}^{n-1}\right)
\tilde f_j \sum_{k=1}^{n}\lambda_k p_r(x_k)\widetilde q_j(x_k)\\
&=&
\sum_{j=0}^{n-m} \tilde f_j \sum_{k=1}^{n}\lambda_k p_r(x_k)\frac{p_j(x_k)}{j+1}\\
&+&
\sum_{j=n-m+1}^{n -1} \tilde f_j \sum_{k=1}^{n}\lambda_k p_r(x_k)\left(
\frac{m+n-j}{2m(j+1)}p_j(x_k)-
\frac{(j-n+m)}{2m(2n-j+1)}p_{2n-j}(x_k)\right) .
\end{eqnarray*}
Consequently, by using \cite{woula_NA}
\begin{equation}\label{ide-p}
p_{2n-j}(x_k)=- p_j(x_k), \qquad k=1,\ldots, n,\qquad n-m<j<n,
\end{equation}
for any $r=0,\ldots, n-1$, we get
\begin{eqnarray*}
c_{n,r}(\tilde f)&=&\sum_{j=0}^{n-m}\tilde f_j < p_r,p_j>\frac 1 {j+1}\\ &+&
\sum_{j=n-m+1}^{n -1}\tilde f_j < p_r,p_j>\frac{m+n-j}{2m(j+1)}+
\sum_{j=n-m+1}^{n -1}\tilde f_j <p_r,p_{j}>\frac{(j-n+m)}{2m(2n-j+1)}.
\end{eqnarray*}
Hence, the orthogonality relation $<p_r,p_{j}>=\delta_{r,j}$ implies that
\[V_n^m\tilde  f(y)= \sum_{r=0}^{n-m}q_r(y)\tilde f_r\left[\frac{1}{r+1}\right]+
\sum_{r=n-m+1}^{n-1}q_r(y)\tilde f_r\left[\frac 1{2m} \left( \frac{m+n-r}{r+1}+\frac{r-n+m}{2n-r+1}\right)\right],
\]
and the statement follows.
\subsection{Proof of Lemma \ref{lemmaKn}.}
Due to the assumption
$\displaystyle \tilde f(y)=\sum_{\ell=0}^{n-1}\tilde f_\ell\ \tilde q_\ell(y)$,
the core of the proof consists in stating that
\begin{equation}\label{Kqtilde}
K \widetilde q_\ell(x_k)= \alpha_\ell p_{\ell-2}(x_k)+\beta_\ell p_\ell(x_k)+\gamma_\ell p_{\ell+2}(x_k),\quad  \ell=0,1,\dots,n -1.
\end{equation}
This is an immediate consequence of (\ref{qtilde-basis}), (\ref{K0}) and (\ref{K-l}) in the case that $\ell=0,\ldots, (n-m)$. In order to prove (\ref{Kqtilde}) also in the case $(n-m)<\ell<n$, we observe that in this case the previous equations yield
 \begin{eqnarray*}
 K \widetilde q_\ell(x_k)&=&
 \frac{n+m-\ell}{2m(\ell+1)}K p_\ell(x_k)-\frac{\ell-n+m}{2m(2n-\ell+1)}K p_{2 n-\ell}(x_k)\\ &=&
  \frac{n+m-\ell}{8m(\ell+1)}\left[-\frac 1 \ell p_{\ell-2}(x_k)+\left(\frac 1 \ell+\frac 1 {\ell+2}\right)p_\ell(x_k)-\frac 1 {\ell+2}p_{\ell+2}(x_k) \right]\\
&-&\frac{\ell-n+m}{8m(2n-\ell+1)}\left[-\frac 1 {2n-\ell} p_{2n-\ell-2}(x_k)+\left(\frac 1 {2n-\ell}+\frac 1 {2n-\ell+2}\right)p_{2n-\ell}(x_k)\right. \\
&-& \left.\frac 1 {2n-\ell+2}p_{2n-\ell+2}(x_k)\right],
  \end{eqnarray*}
and using (\ref{ide-p}), we get
 \begin{eqnarray} \nonumber K \widetilde q_\ell(x_k)=
 &-&\left[\frac{n+m-\ell}{8m\ell(\ell+1)}+\frac{\ell-n+m}{8m(2n-\ell+1)(2n-\ell+2)}\right]
 p_{\ell-2}(x_k) \\ \label{daprima}
 &+&\left[\frac{n+m-\ell}{4m\ell(\ell+2)}+\frac{\ell-n+m}{4m(2n-\ell)(2n-\ell+2)}\right]p_\ell(x_k)\\ & - &\left[
  \frac{n+m-\ell}{8m(\ell+1)(\ell+2)}+\frac{\ell-n+m}{8m(2n-\ell+1)(2n-\ell)}
 \right]p_{\ell+2}(x_k) \nonumber.\end{eqnarray}
By (\ref{daprima}), the statement follows if $(n-m)<\ell\le (n-3)$ and the same holds if $\ell=(n-2)$ being in this case $p_{\ell+2}(x_k)=p_n(x_k)=0$.

If $\ell=(n-1)$ then (\ref{daprima}) implies that
\begin{eqnarray*} K \widetilde q_{n-1}(x_k)=
 &-&\left[\frac{m+1}{8 m n (n-1)}+\frac{m-1}{8m(n+2)(n+3)}\right]p_{ n-3}(x_k) \\
 & +&\left[\frac{m+1}{4m(n-1)(n+1)}+\frac{m-1}{4m(n+1)(n+3)}\right]p_{n -1}(x_k)\\
 & - &\left[\frac{m+1}{8 m n (n+1)}+\frac{m-1}{8 m(n+2)(n+1)}\right]p_{n+1}(x_k),
 \end{eqnarray*}
and taking into account that  $p_{n-1}(x_k)=-p_{n+1}(x_k)$ (cf. (\ref{ide-p})), we get
 \begin{eqnarray*} K \widetilde q_{n-1}(x_k)=
 &-&\left[\frac{m+1}{8 m n (n-1)}+\frac{m-1}{8m(n+2)(n+3)}\right]p_{ n-3}(x_k) \\
 &+&\left[\frac{(m+1)(3n-1)}{8 m n (n-1)(n+1)}+\frac{(m-1)(3n+7)}{8m(n+1)(n+2)(n+3)}\right]p_{ n-1}(x_k),
 \end{eqnarray*}
 which concludes the proof of (\ref{Kqtilde}).

Finally, by (\ref{Kqtilde}), we get the statement as follows
 \begin{eqnarray*}
 K_n^m\tilde f (y)&=&\sum_{j=0}^{n-1}q_j(y)\left[\sum_{k=1}^{n}\lambda_k p_j(x_k)\sum_{\ell=0}^{n-1} \tilde f_\ell K \widetilde q_\ell(x_k)\right]
\\
 &=&\sum_{j=0}^{n-1}q_j(y)\left[\sum_{k=1}^{n}\lambda_k p_j(x_k)\sum_{\ell=0}^{n-1}\tilde f_\ell \bigg(\alpha_\ell p_{\ell-2}(x_k)+\beta_\ell p_\ell(x_k)+\gamma_\ell p_{\ell+2}(x_k)\bigg)\right]\\
 & =& \sum_{j=0}^{n-1}q_j(y)\sum_{\ell=0}^{n-1}\tilde f_\ell \bigg[\alpha_\ell
 <p_{\ell-2}, p_j>+\beta_\ell<p_{\ell}, p_j>+ \gamma_\ell<p_{\ell+2}, p_j>\bigg]\\
 &=&\sum_{j=0}^{n-1}q_j(y)\bigg[
 \alpha_{j+2}\tilde f_{j+2}+\beta_{j}\tilde f_{j}+\gamma_{j-2}\tilde f_{j-2}\bigg].
\end{eqnarray*}
\subsection{Proof of Lemma \ref{lemmaH_n}.}
The statement can be deduced from (\ref{def-Vhy})--(\ref{def-Hn}) and (\ref{mix-pr-0}), as follows
\begin{eqnarray*}
H_n^m \tilde f(y)&=&\sum_{j=0}^{n-1}q_j(y)\left[\sum_{k=1}^{n}\lambda_k p_j(x_k)\widetilde H_n \tilde f(x_k)\right]
\\
&=&\frac 1 \pi\sum_{j=0}^{n-1}q_j(y)\sum_{k=1}^{n}\lambda_k p_j(x_k)\left[\sum_{i=0}^{ n-1}\sum_{s=1}^{n}\lambda_s p_i(x_s)h(x_s,x_k)\int_{-1}^1\tilde f(x) q_i(x)\varphi(x)dx\right]\\
&=&\frac 1 \pi\sum_{j=0}^{n-1}q_j(y)\left[\sum_{i=0}^{n-1}\sum_{k=1}^{n} \sum_{s=1}^{n}\lambda_k \lambda_s p_j(x_k)p_i(x_s)h(x_s,x_k)\sum_{\ell=0}^{n-1} \tilde f_\ell<\widetilde q_\ell, q_i>\right]\\ &=&
\frac 1 \pi\sum_{j=0}^{n-1}q_j(y)\left[\sum_{i=0}^{n-1}\sum_{k=1}^{n} \sum_{s=1}^{n}\lambda_k \lambda_s p_j(x_k)p_i(x_s)h(x_s,x_k)\tilde f_i<\widetilde q_i, q_i>\right].
\end{eqnarray*}
\subsection{Proof of Theorem \ref{th-kappa}}

Since the entries of the matrix $\mathcal{A}'_n$ are infinitesimal, for any arbitrary $\varepsilon>0$ there exists $\nu_\varepsilon$ such that for any $n>\nu_\varepsilon$
the matrix $\mathcal{A}'_n$ can be represented in the following block form
\[\mathcal{A}'_n =: \begin{pmatrix}A_{1,1} & A_{1,2}\\ A_{2,1} & A_{2,2}\end{pmatrix},\]
with
\begin{eqnarray*}A_{1,1}&=&\begin{pmatrix}
\delta_0 & 0 & \alpha_2 &  &     &  & \\
0 & \delta_1 & 0  &  \alpha_3  & 0 & \\
\gamma_0 & 0 &  \delta_2 & 0 & \ddots &  \\
 &  & \ddots & \ddots &  \ddots & \alpha_{\nu_\varepsilon} \\
 &  &  & \ddots & \ddots & 0\\
 &  &  &  \gamma_{\nu_\varepsilon-2} & 0 &  \delta_{\nu_\varepsilon}\\
\end{pmatrix}\in \RR^{(\nu_\varepsilon+1)\times (\nu_\varepsilon+1)} \\ 
A_{1,2}&=&\begin{pmatrix}
0 & 0 & \ldots & 0 & 0   \\
\vdots &  0 & 0 & 0 &   0  \\
0 &  \vdots &   0 & \vdots & \vdots\\
\alpha_{\nu_\varepsilon+1} & 0   & \vdots &  0  & 0 \\
 0 & \alpha_{\nu_\varepsilon+2} &    0 & \ldots & 0
\end{pmatrix}  \in \RR^{(\nu_\varepsilon+1)\times (n-\nu_\varepsilon-1)} \\
A_{2,1}&=& \begin{pmatrix}
0 & 0 & 0  & \gamma_{\nu_\varepsilon-1}   & 0\\
0 &  \vdots & \vdots &  0 & \gamma_{\nu_\varepsilon} \\
\vdots &  \vdots &   0 & \ldots & 0\\
0 & 0   & 0 & \vdots & 0 &  \\
 0 & 0 &    \ldots & 0 & 0
\end{pmatrix} \in \RR^{(n-\nu_\varepsilon-1)\times (\nu_\varepsilon+1) }\\ 
\end{eqnarray*}
\begin{eqnarray*}
 A_{2,2}&=&\begin{pmatrix}
\delta_{\nu_\varepsilon+1} & 0  & \alpha_{\nu_\varepsilon+3} &   &  &  & \\
0 & \delta_{\nu_\varepsilon+2} & 0 &  \alpha_{\nu_\varepsilon+4}  & &  & \\
\gamma_{\nu_\varepsilon+1} & 0 & \delta_{\nu_\varepsilon+3} & 0 &  & &  \\
   &  &   \ddots & \ddots & \ddots &  & 0\\
 &  & \ddots  & \ddots &  \ddots & & \alpha_{n}\\
 &  &    & \ddots & \ddots &  & 0\\
  & &  &   & \gamma_{n-2} & 0 & \delta_{n-1}\\
\end{pmatrix}\in \RR^{(n-\nu_\varepsilon-1)\times (n-\nu_\varepsilon-1)}
\end{eqnarray*}

where $\delta_k=1+\sigma w_{k}+\beta_k,\ k=0,1,\dots, n-1$, and $\alpha_k, \gamma_k,\beta_k, w_{k}$ are defined by Lemmas \ref{lemmaVn} and
\ref{lemmaKn}.  We note that
 the elements $\alpha_i, \gamma_i,\beta_i, w_{i}$ 
are decreasing in modulus and tend to $0$ as   $n\to \infty$, $|\alpha_i|,|\gamma_i|, \beta_i $  with order  $1/n^2$ and $ w_{i}$ with order $1/n$.
The matrix $\mathcal{A}'_n$ can be split into the sum of two  $2\times 2$ block matrices having  homologous blocks of equal dimensions, i.e.
$$\mathcal{R}_n=\begin{pmatrix}A_{1,1} & \mathbf{O}_{(\nu_\varepsilon+1)\times (n-\nu_\varepsilon-1)}\\
 \mathbf{O}_{(n-\nu_\varepsilon-1)\times (\nu_\varepsilon+1)}& I_{(n-\nu_\varepsilon-1)\times (n-\nu_\varepsilon-1)} \end{pmatrix},$$ $$\mathcal{E}_n=\begin{pmatrix}\mathbf{O}_{(\nu_\varepsilon+1)\times (\nu_\varepsilon+1)} & A_{1,2}\\ A_{2,1} & A_{2,2}-I_{(n-\nu_\varepsilon-1)\times (n-\nu_\varepsilon-1)} \end{pmatrix}. $$
Since for $n$ sufficiently large $\|\mathcal{E}_n\|<\varepsilon$, we have
\begin{equation}\label{inequalityAn}\|\mathcal{R}_n\|-\varepsilon\le \| \mathcal{A}'_n\|\le \|\mathcal{R}_n\|+\varepsilon, \quad \forall \epsilon>0, \quad \forall n>\nu_\epsilon,\end{equation}
where we remark that $\|\mathcal{R}_n\|$ is independent of $n>\nu_\varepsilon$, since increasing the order $n$ of $\mathcal{R}_n$,  the order of the identity block increases, while $A_{1,1}$ remains unchanged.
Now, in order to prove the convergence of $\{\|\mathcal{A}'_n\|\}_n$, it will be sufficient to state that $\{\|\mathcal{A}'_n\|\}_n$ is a Cauchy sequence, i.e.
$$\forall \varepsilon>0 \ \exists \nu_\varepsilon:\ \forall n_1,n_2>\nu_\varepsilon, \ \left|\|\mathcal{A}'_{n_1}\|-\|\mathcal{A}'_{n_2}\|\right|<\varepsilon,$$
and this is a consequence of (\ref{inequalityAn}).

Now we prove that also the sequence $\{\|({\mathcal{A}'}_n)^{-1}\|\}_n$  is a Cauchy sequence. For the sake of brevity, we will omit the dimensions of the  blocks involved into the matrices, since they were  specified above.

We have  $$\mathcal{A}'_n=\mathcal{R}_n+\mathcal{E}_n=\mathcal{R}_n(\mathcal{I}_n+\mathcal{R}_n^{-1}\mathcal{E}_n),$$
where

\begin{equation}\label{inv_Bn}\mathcal{R}_n^{-1}=\begin{pmatrix}A_{1,1}^{-1} & \mathbf{O}\\ \mathbf{O} & \mathcal{I} \end{pmatrix},\end{equation}

\begin{equation}\label{mat_blocchi}\mathcal{R}_n^{-1}\mathcal{E}_n =\begin{pmatrix} \mathbf{O} &  A_{1,1}^{-1}A_{1,2}\\
A_{2,1} & A_{22}- I\end{pmatrix},\end{equation}
\begin{eqnarray*}   A_{1,1}^{-1}A_{1,2}&=&A_{1,1}^{-1}[\alpha_{\nu_\varepsilon+1}\mathbf{e}_{\nu_\varepsilon},
\alpha_{\nu_\varepsilon+2}\mathbf{e}_{\nu_\varepsilon+1}, \underbrace{\mathbf{0}, \ldots, \mathbf{0}}]\\ &&\hspace{3.8cm} n-\nu_\varepsilon-3
\end{eqnarray*}
denoting by $ \mathbf{e}_i\in \RR^{\nu_\varepsilon+1}$  the vectors of the canonical bases.

Now, in order to obtain an explicit expression of  $ A_{1,1}^{-1}A_{1,2}$, we need the last two columns of $A_{1,1}^{-1}$. For semplicity we introduce the vectors $$\mathbf{P}^0:=A_{1,1}^{-1}\mathbf{e}_{\nu_\varepsilon}, \quad \mathbf{P}^1:=A_{1,1}^{-1}\mathbf{e}_{\nu_\varepsilon+1},$$
whose elements are denoted as $\mathbf{P}^0_i, \mathbf{P}^1_i,\quad 0\le i\le \nu_\varepsilon-1,$
and for  fixing the ideas, assume $\nu_\varepsilon$ even.
Thus we start from the  $LU$ factorization of $A_{1,1}$
$$A_{1,1}=LU=\begin{pmatrix}
1 &  &  &    &   &  & \\
0 & 1 &  &   &  & \textrm{\huge 0} & \\
v_0 & 0 &  1 &  &   &  &  \\
 & v_1 & 0 &  1 &  &     &  \\
 &  & \ddots & \ddots &  \ddots &  & \\
 & \textrm{\huge 0} &  &  v_{\nu_\varepsilon-2}  & 0 & 1 &
 \end{pmatrix}
  \begin{pmatrix}
d_0  & 0 & \alpha_{2} &  &  &   &   \\
 & d_1  & 0 & \alpha_{3} &  & \textrm{\huge 0} &   \\
 &  &  d_2 & 0 & \ddots & &   \\
 &  &  & \ddots &  \ddots & \alpha_{\nu_\varepsilon}  \\
 &  & \textrm{\huge 0} &   & d_{\nu_\epsilon-1}&  0 \\
 &  &  &  &   & d_{\nu_\varepsilon}
\end{pmatrix},
$$

where $\ d_0=\delta_0,\ d_1=\delta_1,$ and for $k=2,3,\dots \nu_\varepsilon$,
\begin{eqnarray*}\ v_{k-2}&=&\frac{\gamma_{k-2}}{d_{k-2}},\\
 d_k&=&\delta_k-v_{k-2}\alpha_k.\end{eqnarray*}

By induction  it is possible to prove that $d_k:=d_k^{(n)}\to 1$ as   $n\to \infty$, with  $d_k= 1+\mathcal{O}\left(\frac{|\sigma|}n\right)$ and $v_k:=v_k^{(n)}\to 0$ with order $1/n^2$.

Moreover, it is possible to determine by induction

\begin{eqnarray*}
\mathbf{P}^0_{0}&=&\frac{1-\alpha_2 \mathbf{P}^0_2}{d_0},\\
\mathbf{P}^0_{2i+1}&=&0, \quad  i=0,1,\ldots,\frac{\nu_\varepsilon}2-1,\\
 \mathbf{P}^0_{2i}&=&  \frac{(-1)^i\prod_{j=0}^{i-1}v_{2j}-\alpha_{2i+2}\mathbf{P}^0_{2i+2}}{d_{2i}}            \quad  i=\frac{\nu_\varepsilon}2-1, \frac{\nu_\varepsilon}2-2,\ldots,1,\\
  \mathbf{P}^0_{\nu_\varepsilon}&=& \frac{(-1)^\frac{\nu_\varepsilon}2\prod_{j=0}^{\frac{\nu_\varepsilon}2-1}v_{2j}}{d_{\nu_\varepsilon}},
\\ \\
  \mathbf{P}^1_{1}&=&  \frac{1-\alpha_{3}\mathbf{P}^1_{3}}{d_{1}},\\
  \mathbf{P}^1_{2i+1}&=&  \frac{(-1)^i\prod_{j=0}^{i-1}v_{2j+1}-\alpha_{2i+3}\mathbf{P}^1_{2i+3}}{d_{2i+1}}            \quad  i=\frac{\nu_\varepsilon}2-3, \ldots,2,
  \\ \quad \mathbf{P}^1_{2i}&=&0, \quad  i=0,1,\ldots,\frac{\nu_\varepsilon}2,\\
   \mathbf{P}^1_{\nu_\varepsilon-1}&=& \frac{(-1)^{\frac{\nu_\varepsilon}2-1}\prod_{j=0}^{\frac{\nu_\varepsilon}2-2}v_{2j+1}}{d_{\nu_\varepsilon-1}},
\end{eqnarray*}
and taking into account the orders of convergence to $0$  of the sequences $\alpha_k, \delta_k,\gamma_k, v_k$ and the boundedness of $\left\{1/ d_k^{(n)}\right\}_k$ as $n\to \infty$, we can conclude that  under suitable choice  of $\nu_\varepsilon$ the entries of $\mathcal{R}_n^{-1}\mathcal{E}_n$ are so small that
$\|\mathcal{R}_n^{-1}\mathcal{E}_n\|<\varepsilon$, i.e. for any $\varepsilon>0$ we can find $\nu_\varepsilon$ such that for every $n>\nu_\varepsilon$
 $$\mathcal{A}'_n=\mathcal{R}_n(\mathcal{I}_n+\mathcal{R}_n^{-1}\mathcal{E}_n),\qquad \|\mathcal{R}_n^{-1}\mathcal{E}_n\|<\varepsilon.$$
 Taking $\varepsilon<1$ and recalling that $\|(\mathcal{I}_n+\mathcal{R}_n^{-1}\mathcal{E}_n)^{-1}\|\le \frac 1 {1-\|\mathcal{R}_n^{-1}\mathcal{E}_n\|}$
 (see e.g. \cite[Lemma 2.3.3, p. 59]{golub}) we have
\begin{equation}\label{eq_bound}\frac{\|\mathcal{R}_n^{-1}\|}{1+\varepsilon}\le \|{(\mathcal{A}'}_n)^{-1}\|\le \frac{\|\mathcal{R}_n^{-1}\|}{1-\varepsilon},\quad \forall n>\nu_\varepsilon.\end{equation}
As we have said before  $\|\mathcal{R}_n\|$ is independent of $n>\nu_\varepsilon$ as well as,   in view of  (\ref{inv_Bn}),  $\|\mathcal{R}_n^{-1}\|$ is too.
 Therefore $ \|({\mathcal{A}'_n})^{-1}\|$ is  bounded  $\forall n\in\NN$, i.e. there exists a positive constant $M$ s.t. \begin{equation}\label{boundA}\|({\mathcal{A}'}_n)^{-1}\|<M,\quad  \forall n\in \NN, \quad M\neq M(n)\end{equation} and as consequence, by the left hand side bound in (\ref{eq_bound})
$$\|\mathcal{R}_n^{-1}\|\le \|({\mathcal{A}'}_n)^{-1}\|(1+\varepsilon)<2M,\quad  \forall n>\nu_\varepsilon.$$
By  (\ref{boundA}) and taking into account last inequality,  we get for any $n_1,n_2>\nu_\varepsilon$ and $0<\varepsilon<1$
$$ \left|\|({\mathcal{A}'}_{n_1})^{-1}\|-\|({\mathcal{A}'}_{n_2})^{-1}\|\right|\le 2\varepsilon\frac{\|\mathcal{R}_n^{-1}\|}{1-\varepsilon^2}\le
 \frac{4M\varepsilon}{1-\varepsilon^2}$$
 by which
 $$\lim_{n_1,n_2\to\infty}\left|\|({\mathcal{A}'}_{n_1})^{-1}\|-\|({\mathcal{A}'}_{n_2})^{-1}\|\right|=0.$$
 Since also $\{\|({\mathcal{A}'}_n)^{-1}\|\}_n$ is a Cauchy sequence, the thesis follows taking into account the convergence of  the sequence
$\{\|\mathcal{A}'_n\|\}_n$.

{\small{\it Maria Carmela De Bonis, {Department of Mathematics, Computer Science and Economics, University of Basilicata, Via dell'Ateneo Lucano 10, 85100 Potenza, Italy.}\\  mariacarmela.debonis@unibas.it.

\vspace{0.5cm}

Donatella Occorsio {Department of Mathematics, Computer Science and Economics, University of Basilicata, Via dell'Ateneo
Lucano 10, 85100 Potenza, Italy. \\ donatella.occorsio@unibas.it.}

\vspace{0.5cm}

Woula Themistoclakis, C.N.R. National Research Council of Italy, IAC Institute for Applied Computing ``Mauro Picone'',\\ Via P. Castellino, 111, 80131 Napoli, Italy. \\ {woula.themistoclakis@cnr.it.}

}}
\end{document}